\documentclass[11pt]{article}
\usepackage{latexsym}
\usepackage{epsfig}
\usepackage{amssymb, amsmath, multicol}
\usepackage[dvips]{color}
\usepackage{graphicx}
\renewcommand{\baselinestretch}{1.5} 

\makeatletter

\def\singlespace{\def\baselinestretch{1.3}\@normalsize}

\newtheorem{lemma}{Lemma}
\newtheorem{proposition}{Proposition}
\newtheorem{theorem}{Theorem}

\renewcommand{\theequation}{\arabic{equation}}
\renewcommand{\thefootnote}{\fnsymbol{footnote}}
\makeatother

\newcounter{assump}

\newcommand{\bb}{\mbox{\bf b}}
\newcommand{\ba}{\mbox{\bf a}}
\newcommand{\bA}{\mbox{\bf A}}

\newcommand{\bt}{\mbox{\bf t}}

\newcommand{\bW}{\mbox{\bf W}}

\newcommand{\bSig}{\mbox{\boldmath $\Sigma$}}
\newcommand{\bc}{\mbox{\bf c}}
\newcommand{\be}{\mbox{\bf e}}
\newcommand{\vecc}{\text{vec}}
\newcommand{\vech}{\text{vech}}
\newcommand{\APE}{\mbox{APE}}

\def\toD{\overset{\text{D}}{\longrightarrow}}
\def\toP{\overset{\text{P}}{\longrightarrow}}
\def\toas{\overset{\text{a.s.}}{\longrightarrow}}

\newcommand{\bvarepsilon}{\mbox{\boldmath $\varepsilon$}}

\newcommand{\bY}{\mbox{\boldmath $Y$}}
\newcommand{\bX}{\mbox{\boldmath $X$}}

\newcommand{\bM}{\mbox{\boldmath $M$}}

\newcommand{\bmu}{\mbox{\boldmath $\mu$}}
\newcommand{\bsigma}{\mbox{\boldmath $\sigma$}}
\newcommand{\PE}{\mbox{PE}}

\newcommand{\var}{\mbox{var}}

\newcommand{\cov}{\mbox{cov}}

\def\toD{\overset{\text{D}}{\longrightarrow}}
\def\toP{\overset{\text{P}}{\longrightarrow}}
\newcommand{\etal}{{\it et al.}}

\begin{document}

\title{\bf Aggregation of Nonparametric Estimators for Volatility Matrix
\thanks{Financial support
from the NSF under grant DMS-0532370 is gratefully acknowledged. We
are grateful to Lars Hansen, Robert Kimmel, and conference
participants of Financial Mathematics Workshop at SAMSI; the 2006
Xiamen Financial Engineering and Risk Management Workshop;  the 69th
IMS Annual Meeting for helpful comments.}
\date{November 1, 2006}
\vspace{0.4 in}
\author{Jianqing Fan\thanks{Princeton, NJ
08544. Phone: (609) 258-7924. E-mail: jqfan@princeton.edu.}\\
{\small Department of Operations Research and Financial
Engineering}\\{\small Princeton University} \and Yingying Fan
\thanks{Corresponding author. Princeton, NJ
08544. Phone: (609) 258-7383. E-mail:
yingying@princeton.edu.}\\{\small Department of Operations Research
and Financial Engineering}\\{\small Princeton University} \and
Jinchi Lv \thanks{Princeton, NJ 08544. Phone: (609) 258-9433.
E-mail: jlv@princeton.edu.}\\
{\small Department of Mathematics}\\{\small Princeton University}} }

\maketitle

\vspace{0.25 in}

\noindent Address correspondence to Yingying Fan, Department of
ORFE, Princeton University, Princeton, NJ 08544. Email:
yingying@princeton.edu.

\newpage
\noindent{\bf ABSTRACT}\\
\noindent An aggregated method of nonparametric estimators based on
time-domain and state-domain estimators is proposed and studied. To
attenuate the \emph{curse of dimensionality}, we propose a factor
modeling strategy. We first investigate the asymptotic behavior of
nonparametric estimators of the volatility matrix in the time domain
and in the state domain. Asymptotic normality is separately
established for nonparametric estimators in the time domain and
state domain. These two estimators are asymptotically independent.
Hence, they can be combined, through a dynamic weighting scheme, to
improve the efficiency of volatility matrix estimation. The optimal
dynamic weights are derived, and it is shown that the aggregated
estimator uniformly dominates volatility matrix estimators using
time-domain or state-domain smoothing alone. A simulation study,
based on an essentially affine model for the term structure, is
conducted, and it demonstrates convincingly that the newly proposed
procedure outperforms both time- and state-domain estimators.
Empirical studies further endorse the advantages of our aggregated
method.

\bigskip


\bigskip

\noindent Covariance matrices are fundamental for risk management,
asset pricing, proprietary trading, and portfolio managements.
In forecasting a future event such as the volatility matrix, two
pieces of information are frequently consulted. Based on the recent
history, one uses a form of local average, such as the moving
average, to predict the volatility matrix. This approach localizes
in time and uses the smoothness of the volatility matrix as a
function of time. It ignores completely the historical information,
which is related to the current prediction through a
stationarity{\renewcommand{\thefootnote}{\arabic{footnote}}\footnotemark\footnotetext{By
``stationarity" we do not mean that the process is strongly
stationary, but has some structural invariability over time. For
example, the conditional moment functions do not vary over time.}
assumption. On the other hand, one can predict a future event by
consulting the historical information with similar scenarios. This
approach basically localizes in the state variable and depends on
the stationarity assumption. For example, by localizing on a few key
financial factors, one can compute the volatility matrix using the
historical information.  This results in a nonparametric estimate of
the volatility matrix using state-domain smoothing. See, for
example, Anderson, Bollerslev and Diebold (2002) for a unified
framework of interpreting both parametric and nonparametric
approaches for volatility measurement.

The aforementioned two estimators are weakly correlated, as they use
data that are quite far apart in time. They can be combined to
improve the efficiency of the volatility matrix estimation.  This
results in an aggregated estimator of the volatility matrix.
Three challenges arise in the endeavor: the \emph{curse of
dimensionality}, the choice of dynamic weights, and the mathematical
complexity.

Due to the curse of dimensionality, surface smoothing techniques are
not very useful in practice when there are more than two or three
predictor variables. An efficient dimensionality reduction process
should be imposed in state-domain estimation. An introduction to
some of these approaches, such as \emph{additive modeling, partially
linear modeling, modeling with interactions}, and \emph{multiple
index models}, can be found in Fan and Yao (2003).

In this paper, we propose a \emph{factor modeling} strategy to
reduce the dimensionality 
in the state domain smoothing. Specifically, to estimate the
covariance matrix among several assets, we first find a few factors
that capture the main price dynamics of the underlying assets.
Regarding the covariance matrix as a smooth function of these
factors, the covariance matrix can be computed via localizing on the
factors.

\begin{center}
Figure 1 here.
\end{center}

Our approach is particularly appealing for the yields of bonds, as
they are often highly correlated, which makes the choice of the
factors relatively easy. 
To elucidate our idea, consider the weekly data on the yields of
treasury bills and bonds with maturities 1 year, 5 years, and 10
years presented in Figure 1. We choose the 5-year yield process as
the single factor. Suppose that the current time is January 14, 2000
and the current interest rate of the 5-year treasury bond is 6.67\%,
corresponding to time index $t = 1986$. One may estimate the
volatility matrix based on the weighted squared differences in the
past 104 weeks. This corresponds to time-domain smoothing, using the
small vertical stretch of data shown in Figure 1(a). On the other
hand, one may also estimate the volatility matrix using the
historical data with interest rates approximately 6.67\%, say,
$6.67\% \pm .20\%$. This corresponds to localizing in state domain
and is indicated by the horizontal bar in Figure 1(a). Figures 1(b)
and 1(c) present scatter plots of the yield differences
$X^{1\text{yr}}_t - X^{1\text{yr}}_{t-1}$ for the 1-year bill
against the yield differences
$X^{10\text{yr}}_t-X^{10\text{yr}}_{t-1}$ for the 10-year bond,
using respectively the data localizing in the time and state
domains. The associated regression lines of the time- and
state-domain data are also presented. The scatter plots give two
estimates of the conditional correlation and conditional variance of
the volatility matrix for the week of $t=1986$. They are weakly
dependent as the two scatter plots use data that are quite far apart
in time.

Let $\widehat{\bSig}_{T,t}$ and $\widehat{\bSig}_{S,t}$ be the
estimated volatility matrices based on data localizing in the time
and state domains, respectively. For example, they can be the sample
covariance matrices for the data presented in Figures 1(b) and 1(c),
respectively for $t=1986$. To fully utilize these two estimators, we
introduce a weight $w_t$ and define an aggregated
estimator\footnote{Ledoit and  Wolf (2003) introduce a shrinkage
estimator by combining the sample covariance estimator with that
derived from the CAPM.  Their procedure intends to improve estimated
covariance matrix by pulling the sample covariance towards the
estimate based on the CAPM.  Their basic assumption is that the
return vectors are i.i.d. across time.  This usually holds
approximately when the data are localized in time.  In this sense,
their estimator can be regarded as a time-domain estimator.} as
\begin{align}\label{a1}
\widehat{\bSig}_{A,t}=\omega_t\widehat{\bSig}_{S,t} +
(1-\omega_t)\widehat{\bSig}_{T,t}.
\end{align}
The weight function $\omega_t$ is always between 0 and 1, and it can
be an adaptive random process which is observable at time $t$. Due
to the weak dependence between the original two estimators, the
aggregated estimator is always more efficient than either of the
time- and state-domain estimators.

An interesting question is the choice of the dynamic weight
$\omega_t$. Suppose we have a portfolio with allocation vector
$\ba$. Then the aggregation method gives us the following estimate
of the portfolio variance:
\begin{align}\label{a7}
\ba^T\widehat{\bSig}_{A,t}\ba=\omega_t\ba^T\widehat{\bSig}_{S,t}\ba
+ (1-\omega_t)\ba^T\widehat{\bSig}_{T,t}\ba.
\end{align}
Since $\widehat{\bSig}_{S,t}$ and $\widehat{\bSig}_{T,t}$ are
asymptotically
independent\renewcommand{\thefootnote}{\arabic{footnote}}\footnotemark\footnotetext{
We prove in Section 4 that $\widehat{\bSig}_{S,t}$ and
$\widehat{\bSig}_{T,t}$ are asymptotically independent, and thus
they are close to be independent in finite sample. In the following,
by ``nearly independent'' and ``almost uncorrelated'', we mean the
same.}, the optimal weight in terms of minimizing the variance of
$\ba^T\widehat{\bSig}_{A,t}\ba$ is
\begin{align}\label{a6}
\omega_{\text{opt},t}=\frac{\var(\ba^T\widehat{\bSig}_{T,t}\ba)}
{\var(\ba^T\widehat{\bSig}_{S,t}\ba)+\var(\ba^T\widehat{\bSig}_{T,t}\ba)}.
\end{align}
Indeed, our asymptotic result in Section 3 shows that the optimal
weight admits a simple and explicit form, independent of $\ba$. This
makes our implementation very easy.

The above approach is data analytic in the sense that it is always
operational. To appreciate our idea, we will introduce a
mathematical model for the data-generating process in Section 1. And
then in the following several sections we formally show that the
aggregated estimator has the desired statistical properties.

\section{Model and Assumptions}

Let $\bW_t =(W_1^t, \cdots, W_t^m)^T$ and $\bW=\{\bW_t,\
\mathcal{F}_t^W;\ 0\leq t<\infty\}$ be an $m$-dimensional standard
Brownian motion. Consider the following $d$-dimensional diffusion
process
\begin{equation}\label{a2}
d\bX_t=\bmu_tdt+\bsigma_td\bW_t,
\end{equation}
where $\bX_t= (X_t^1,\cdots, X_t^d)^T$,  $\bmu_t$ is a $d\times1$
predictable vector process, and $\bsigma_t$ is a $d\times m$
predictable matrix process, depending only on $\bX_t$. Here, $m$ can
be different from $d$. This is a widely used model for asset prices
and the yields of bonds. This family of models includes famous ones
such as multivariate generalizations of both Vasicek (1977) and Cox,
Ingersoll and Ross (1985).

Under model (\ref{a2}), the diffusion matrix is $\bSig_t=
\bsigma_t\bsigma_t^T$. As mentioned before, when $d \geq 2$, the
so-called \emph{curse of dimensionality} makes implementation
hard. 
To reduce the dimensionality, we introduce a scalar
factor $f_t$ and model the drift and diffusion processes as
$\bmu_t=\bmu(f_t)$ and $\bsigma_t=\bsigma(f_t)$, where
$\bmu(\cdot)=\{\mu_i(\cdot)\}_{1\leq i\leq d}$ is a $d\times1$ Borel
measurable vector and $\bsigma(\cdot)=\{\sigma_{ij}(\cdot)\}_{1\leq
i\leq d,1\leq j\leq m}$ is a $d\times m$ Borel measurable matrix.
Then model (\ref{a2}) becomes
\begin{equation}\label{a3}
dX^i_t=\mu_i(f_t)dt+\sum_{j=1}^m\sigma_{ij}(f_t)dW_t^j,\quad 1\leq
i\leq d .
\end{equation}
In this model, the diffusion matrix is
$\bSig(f_t)=\bsigma(f_t)\bsigma(f_t)^T$. See also Engle, Ng and
Rothchild (1990) for a similar strategy.

We introduce some stochastic structure on $f_t$ by assuming that
$f_t$ is the solution to the following stochastic differential
equation (SDE):
\begin{equation}\label{a4}
df_t=a(f_t)dt+\sum_{j=1}^mb_j(f_t)dW_t^j,
\end{equation}
where $a(\cdot)$ and $b_1(\cdot),\ b_2(\cdot),\cdots, \ b_m(\cdot)$
are unknown functions. In some situations like modeling bond
yields\footnote{In practice, one can take the yields process with
median term of maturity as the driving factor, as this bond is
highly correlated to both short-term and long-term bonds.}, the
factor $f_t$ can be chosen as one of the bond yields, i.e., $f_t$ is
one of the coordinates of $\bX_t$. But in general, $f_t$ may be
different from any coordinate of $\bX_t$, and the theoretical
studies in this paper apply to both cases. The data are observed at
times $t_i = t_0 + i \Delta$, $i=0,1,\cdots,N$, with sampling
interval $\Delta$, resulting in vectors $\{\bX_{t_i}, i=0,1, \cdots,
N\}$ and $\{f_{t_i}, i=0,1, \cdots, N\}$. This model is reasonable
for the yields of bonds with different maturities since they are
highly correlated. Thus, localizing on all the yields processes in
the state domain results in approximately the same data set as
localizing on only one of the yields
processes. 
In addition, our study can be generalized to the multi-factor case
without much extra difficulty. We will focus on the one-factor
setting for simplicity of presentation.


Let $\bY_i=(\bX_{t_{i+1}}-\bX_{t_i})\Delta^{-1/2}$, and denote by
$Y_i^1,\ Y_i^2,\ \cdots,\ Y_i^d$ the coordinates of $\bY_i$. Then,
by the Euler scheme, we have
\begin{equation}\label{a5}
\bY_i\approx\bmu(f_{t_i})\sqrt{\Delta}+\bsigma(f_{t_i})\bvarepsilon_{t_i},
\end{equation}
where $\bvarepsilon_{t_i}$ follows the $m$-dimensional standard
Gaussian distribution. The conditional covariance matrix of $\bX$ at
time $t_i$ can be approximated by $\Delta\bSig(f_{t_i})$ (see Fan
and Zhang, 2003). Hence, the estimate of the conditional covariance
matrix is almost equivalent to the estimate of the diffusion matrix
$\bSig(\cdot)$. Fan and Zhang (2003) study the impact of the order
of difference on nonparametric estimation.  They found that while
higher order can possibly reduce approximation errors, it increases
variances of data substantially. They recommended the Euler scheme
(\ref{a5}) for most practical situations.

To use time-domain information, it is necessary to assume that the
sampling frequency $\Delta$ converges to zero so that the biases in
time-domain approximations are negligible. As a result, we face the
challenge of developing asymptotic theory for the diffusion model
(\ref{a3}). Both nonparametric estimators in the time domain and
state domain need to be investigated. Pioneering efforts on
nonparametric estimation of drift and diffusion include Jacod
(1997), Jiang and Knight (1997), Arfi (1998),
Gobet (2002), Bandi and Philips (2003), Cai and Hong (2003), Bandi
and Moloche (2004), and Chen and Gao(2004). Arapis and Gao (2004)
investigate the mean aggregated square errors of several methods for
estimating the drift and diffusion, and compare their performances.
A\"{\i}t-Sahalia and Mykland (2003, 2004) study the effects of
random and discrete sampling when estimating continuous-time
diffusions. Bandi and Nguyen (1999) investigate small sample
behaviors of nonparametric diffusion parameters. See Bandi and
Phillips (2002) for a survey of recently introduced techniques for
identifying nonstationary continuous-time processes. As long as the
time horizon is long, the diffusion matrix can be estimated with low
frequency data (say, finite $\Delta^{-1}$). See, for example, Hansen
{\it et al.} (1998) for the spectral method, Kessler and S{\o}rensen
(1999) for parametric models, and Gobet {\it et al.} (2004) for
specific univariate nonparametric diffusions.

To facilitate our future presentation, we make the following
assumptions:


\noindent{\bf\em Assumption 1}. (Global Lipschitz and linear growth
conditions). There exists a constant $k_0\ge 0$ such that
\begin{equation}\label{pf1}
\|\bmu(x)-\bmu(y)\|+\|\bsigma(x)-\bsigma(y)\|\leq k_0 |x-y|,
\end{equation}
\begin{equation*}
\|\bmu(x)\|^2+\|\bsigma(x)\|^2\leq k_0^2(1+x^2),
\end{equation*}
for any $x, y\in \mathbb{R}$. 
Also, with $\bb(\cdot)=(b_1(\cdot),\ b_2(\cdot),\ \cdots,\
b_m(\cdot))^T$, assume that
$$
|a(x)-a(y)|+\|\bb(x)-\bb(y)\|\leq k_0 |x-y|.
$$

\bigskip

\noindent{\bf\em Assumption 2}. Given any time point $t>0$, there
exists a constant $L>0$ such that $E|\mu_i(r_s)|^{4(q_0+\delta)}\leq
L$ and $E|\sigma_{ij}(r_s)|^{4(q_0+\delta)}\leq L$ for any $s\in
[t-\eta,t]$ and $1\leq i,j\leq d$, where $\eta$ is some positive
constant, $q_0$ is an integer not less than $1$, and $\delta$ is
some small positive constant.


\bigskip

\noindent{\bf\em Assumption 3}. The solution $\{f_{t} \}$ of model
(\ref{a4}) is a stationary Markov process and real ergodic. For
$t\geq0$, define the transition operator by:
$$
(H_t g)(a)=E(g(f_t)|f_0=a),\ a\in R,
$$
where $g(\cdot)$ is any Borel measurable bounded function on
$\mathbb{R}$. Suppose $H_t$ satisfy the $G_2$ condition of
Rosenblatt (1970), i.e., there is some $s>0$ such that
$$
|H_s|_2=\sup_{\{g,\
Eg(X)=0\}}\frac{E^{1/2}(H_sg)^2(X)}{E^{1/2}g^2(X)}\leq \alpha<1.
$$

\bigskip

\noindent{\bf\em Assumption 4}. The conditional density
$p_\ell(y|x)$ of $f_{t_{i+\ell}}$ given $f_{t_i}$ is continuous in
the arguments $(y,x)$ and is bounded by a constant independent of
$\ell$. The time-invariant density function  $p(x)$ of the process
$f_t$ is bounded and continuous.

\bigskip

\noindent{\bf\em Assumption 5}. The kernel $K(\cdot)$ is a
continuously differentiable, symmetric probability density function
satisfying
\begin{align}\label{pf2}
\int |x^jK'(x)|dx<\infty,\ j=0,1,\cdots,5,
\end{align}
\begin{align}\label{044}
\mu_i=\int x^iK(x)dx<\infty,\ i=0,1,\cdots,4,
\end{align}
and
$$
\nu_0=\int K^2(x)dx<\infty.
$$

Let $\{\mathcal{F}_t\}$ be the augmented filtration defined in Lemma
2 of Appendix. Assumption 1 ensures that there exist continuous,
adapted processes $\bX=\{\bX_t,\in \mathcal{F}_t;\ 0\leq t<\infty\}$
and $f=\{f_t\in \mathcal{F}_t;\ 0\leq t<\infty\}$, which are strong
solutions to SDEs (\ref{a2}) and (\ref{a4}) respectively, provided
that the initial values $\bX_0$ and $f_0$ satisfy
$E\|\bX_0\|^2<\infty$ and $E|f_0|^2<\infty$, and  are independent of
Brownian motion $\bW$ (see, e.g., Chapter 5, Theorem 2.9 of Karatzas
and Shreve, 1991). Assumption 2 indicates that, given any time point
$t>0$, there is a time interval $[t-\eta,t]$ on which the drift and
volatility functions have finite $4(q_0+\delta)$-th moments.
Assumption 3 says that $f_t$ is stationary and ergodic and satisfies
some mixing condition (see Fan, 2005), which ensures that $f_t$ is
Harris recurrent. For the stationarity assumption of $f_t$ to be
true, see Hansen and Scheinkman (1995) for conditions. Assumption 4
imposes some constraints on the transition density of $f_t$.
Assumption 5 is a regularity condition on the kernel function. For
example, the
commonly used Gaussian kernel satisfies it. 

With the above theoretical framework and assumptions, we will
formally demonstrate that the nonparametric estimators using the
data localizing in time and in state are asymptotically jointly
normal and independent. This
gives a formal theoretical justification 
and serves as the theoretical foundation for the idea that the
time-domain and state-domain nonparametric estimators can be
combined to yield a more efficient volatility matrix estimator.


\section{DIFFUSION MATRIX ESTIMATION USING RECENT INFORMATION}\label{sec1}

The time-domain method has been extensively studied in the
literature. See, for example, Robinson (1997),
H\"ardle \etal\ (2002), Fan, Jiang, Zhang and Zhou (2003), and
Mercurio and Spokoiny (2004),
among others.  
A popular time-domain method, the moving average estimator is
defined as
\begin{align}\label{b1}
\widehat{\bSig}_{MA,t}=\frac{1}{n}\sum_{i=1}^{n}\bY_{t-i}\bY_{t-i}^T,
\end{align}
where $n$ is the size of the moving window. This estimator ignores
the drift component and utilizes $n$ local data points. An extension
of the moving average estimator is the exponential smoothing
estimator, which is defined as
\begin{equation}\label{b2}
\widehat{\bSig}_{ES,t}=(1-\lambda)\sum_{i=1}^\infty\lambda^{i-1}\bY_{t-i}\bY_{t-i}^T,
\end{equation}
where $\lambda$  is a smoothing parameter controlling the size of
the local neighborhood. RiskMetrics of J.P. Morgan (1996), which is
used for measuring the risks of financial assets, recommends
$\lambda = 0.94$ and $\lambda = 0.97$ when one uses (\ref{b2}) to
forecast the daily and monthly volatility, respectively.

The exponential smoothing estimator (\ref{b2}) is one type of
rolling sample variance estimator. See Foster and Nelson (1996) for
more information about rolling sample variance estimators. Estimator
(\ref{b2}) is also related to the multivariate GARCH model in the
literature. Note that when $\Delta$ is very small, the first term on
the right hand side of (\ref{a5}) can be ignored. Thus (\ref{a5})
and (\ref{b2}) can be written as
$$
\bY_i \approx \bsigma(f_{t_i}) \bvarepsilon_i,
$$
$$
\bSig_{t_i}=(1-\lambda)\bY_{i-1}\bY_{i-1}^T+\lambda\bSig_{t_{i-1}},
$$
where $\bSig_{t_i}=\bsigma(f_{t_i})\bsigma(f_{t_i})^T$, which
reminisces the IGARCH model.


The exponential smoothing estimator in (\ref{b2}) is a weighted sum
of squared returns prior to time $t$. Since the weight decays
exponentially, it essentially uses recent data. To explicitly
account for this, we use a slightly modified version:
\begin{equation}
\widehat{\bSig}_{T,t}=\frac{1-\lambda}{1-\lambda^n}\sum_{i=1}^n
\lambda^{i-1}\bY_{t-i}\bY_{t-i}^T.   \label{b3}
\end{equation}
Here, as in the case of the moving average estimator in (\ref{b1}),
$n$ is a smoothing parameter controlling explicitly the window
width, and $\lambda$ acts like a kernel weight which may depend on
$n$. For example, when $\lambda=1-\frac{\tau}{n}$ with $\tau$ a
positive constant, besides the normalization factor
$\frac{1-\lambda}{1-\lambda^n}$, the first data point $Y_{t-1}$
receives weight 1, while the last point $Y_{t-n}$ receives
approximately  weight $e^{-\tau}$. In particular, when $\lambda=1$,
it becomes the moving average estimator (\ref{b1}).

Before going into the details, we first introduce some notations and
definitions. Let $A=(a_{ij})$ be an $m\times n$ matrix. By
$\vecc(A)$ we mean the $mn\times1$ vector formed by stacking the
columns of $A$. 
If $A$ is also symmetric, we vectorize the lower half of $A$ and
denote the vector by $\vech(A)$. These notations are consistent with
Bandi and Moloche (2004).
It is not difficult to verify that there exists a unique $m^2\times
m(m+1)/2$ matrix $D$ with elements 0 and 1, such that
\begin{align*}
P_D\text{vec}(A)=\text{vech}(A),
\end{align*}
where $P_D=(D^TD)^{-1}D^T$. Another useful definition is the
Kronecker product of two matrices $A$ and $B$, which is defined as
$A\otimes B=(a_{ij}B)$.

Since the estimator $\widehat{\bSig}_{T,t}$ is symmetric, we only
need to consider the asymptotic normality of the linear combination
of the vector $\vech(\widehat{\bSig}_{T,t})$:
\begin{align}\label{b4}
\widehat{U}_{T,t}\equiv\bc^T\text{vech}\widehat{\bSig}_{T,t}
=\frac{1-\lambda}{1-\lambda^n}
\sum_{i=1}^n\lambda^{i-1}\sum_{k=1}^d\sum_{\ell=1}^kc_{k\ell}Y_{t-i}^kY_{t-i}^\ell,
\end{align}
where $\bc=(c_{1,1},\ c_{2,1},\ c_{2,2},\ c_{3,1},\ \cdots,\
c_{d,d})^T$ is a constant vector. 
%
%

\medskip

\begin{proposition}\label{P1} Under
Assumptions 1 and 2, for almost every sample path, we have
\begin{equation}\label{b5}
\|\bsigma(r_s)-\bsigma(r_u)\|\leq K|s-u|^q,\quad s,u\in[t-\eta,t],
\end{equation}
where $q=(2q_0-1)/(4q_0)$, $q_0$ is the integer in Assumption 2, and
the coefficient $K$ satisfies $E[K^{4(q_0+\delta)}]<\infty$ with
$\delta$ a positive constant.
\end{proposition}

\noindent{\bf\em Remark 1}.\quad Proposition 1 shows the continuity
of $\bsigma(r_s)$ as a function of time $s$, which is the foundation
of time-domain estimation. In the proof of Proposition \ref{P1}, we
only used Assumption 2 and the condition
$\|\bsigma(x)-\bsigma(y)\|\leq k_0|x-y|$ with $k_0$ a positive
constant. Assumption 1 is made to ensure the existence of a solution
to model (\ref{a3}).


\begin{theorem} \label{T1}  Suppose that $n\rightarrow\infty,\
n\Delta^{2q/(2q+1)}\rightarrow0$, and Assumptions 1 and 2 hold at
time $t$. If the limit $\tau=\lim\limits_{n\to \infty}n(1-\lambda)$
exists, then 
given $f_t = x$, the conditional distribution of
$\mathrm{vech}(\widehat{\bSig}_{T,t})$ is asymptotically normal,
i.e.,
\begin{align*}
\sqrt{n}\ \mathrm{vech}(\widehat{\bSig}_{T,t}-\bSig(x)) \toD
N\left(0,\frac{\tau(1+e^\tau)}{(e^\tau-1)}\Lambda(x)\right),
\end{align*}
where $\Lambda(x)=P_D^T\{\bSig(x)\otimes\bSig(x)\}P_D$.
\end{theorem}

\medskip

Note that all data used in the estimator (\ref{b3}) is within
$n\Delta$ away from time $t$.  According to Proposition 1, the
approximation error of (\ref{b3}) is at most of order $O
((n\Delta)^q)$, which together with the condition $n
\Delta^{2q/(2q+1)} \to 0$ in Theorem \ref{T1} guarantees that the
bias is of order $o(n^{-1/2})$.

\bigskip

\section{DIFFUSION MATRIX ESTIMATION USING HISTORICAL INFORMATION}\label{sec2}

The diffusion matrix in (\ref{a2}) can also be regarded as a
nonparametric regression given $f_t = x$. See for example its first
order approximation (\ref{a5}). Therefore, it can be estimated by
using the historical information via localizing on the state
variable $f_t$, as illustrated in Figure 1. The local linear
smoother studied in Stanton (1997) will be employed. This technique
has several nice properties, such as asymptotic minimax efficiency
and design adaptation. Further, it automatically corrects edge
effects and facilitates bandwidth selection (Fan and Yao, 2003).

In the construction of the state-domain estimator, we will use the
$N-1$ data points right before the current time $t$, i.e., the
historical data $\{(f_{t_i},\ \bY_i),\ i=0,\ 1,\ \cdots,\ N-1\}$.

It can be shown that the diffusion matrix has the standard
interpretation in terms of infinitesimal conditional moments, that
is,
\begin{align*}
E[Y^i_kY^j_k|f_{t_k}=x_0]= v_{ij}(x_0)+O(\Delta).
\end{align*}

For a given kernel function\footnote{The kernel function is a
probability density, and the bandwidth is its associated scale
parameter.  Both of them are used to localize the linear regression
around the given point $x_0$.  The commonly used kernel functions
are the Gaussian density and the Epanechnikov kernel $K(x) = 0.75
(1-x^2)_+$.} $K$ and a bandwidth $h$, the local linear estimator
$\hat{\beta}^{ij}_0$ of $v_{ij}(x_0)$ is obtained by minimizing the
objective function
\begin{align}\label{c1}
\sum_{k=0}^{N-1}\{Y^i_kY^j_k+\beta^{ij}_0+(f_{t_k}-x_0)\beta^{ij}_1
\}K_h(f_{t_k}-x_0)
\end{align}
over $\beta^{ij}_0$ and $\beta^{ij}_1$. Let
\begin{align}\label{c2}
W_\ell(x) = \sum_{k=0}^{N-1}(f_{t_k} -x)^\ell K_h(f_{t_k}-x)
\end{align}
and
\begin{align}\label{c3}
\quad w_k(x)= K_h(f_{t_k}-x) \{W_2(x) -
        (f_{t_k}-x)W_1(x) \}/\{ W_0(x) W_2(x) - W_1(x)^2\}.
\end{align}
Then the local linear estimator in (\ref{c1}) can be expressed as
\begin{equation}\label{c4}
\widehat{\bSig}_{S,t}(x) = \sum_{k=0}^{N-1} w_k(x) \bY_k\bY_k^T.
\end{equation}
This estimator depends only on the historical data (horizontal bar
in Figure 1), and relies on the structure invariability.

The above weight function $w_k(x)$ is called an ``equivalent kernel"
in Fan and Yao (2003). Expression (\ref{c4}) reveals that the
estimator $\widehat{\bSig}_{S,t}(x)$ is very much like a
conventional kernel estimator except that the ``kernel" $w_k(x)$
depends on the design points and locations.

Before establishing the asymptotic normality of
$\widehat{\bSig}_{S,t}(x)$, we first investigate the asymptotic
property of $W_\ell(x)$.

\begin{proposition} \label{P2}
Suppose $\Delta\rightarrow0$, $N\Delta\rightarrow\infty$, and
$\frac{1}{h}\sqrt{\Delta\log\Delta^{-1}}=o(1)$. Under Assumptions
3--5, we have
\begin{align}\label{c5}
 W_\ell(x)=Nh^\ell\{p(x)\mu_\ell+o_{a.s.}(1)\},\ \ell=0,1,2,3.
\end{align}
\end{proposition}

\medskip

The results of Proposition \ref{P2} are similar to those in Section
6.3.3 of Fan and Yao (2003, p.237), but the proofs are completely
different, as we have a highly correlated sample $\{f_{t_i}\}$ here.
The high correlation makes their proof fail in our case. To attack
this problem, we invoke the local time. The definition and some
preliminary results of local time can be found in Revuz and Yor
(1999, p.221). For the multifactor situation, the local time
generally does not exist. However, by using the occupation time of
Bandi and Moloche (2004), our results can be generalized to the
multifactor situation.


\begin{theorem}\label{T2} Suppose $\Delta\rightarrow0$, $N\Delta\rightarrow\infty$,
$h=O(N^{-1/5})$, and $\frac{1}{h}\sqrt{\Delta\log\Delta^{-1}}=o(1)$.
Moreover, suppose that $\bSig(\cdot)$ is twice differentiable. Under
Assumptions\renewcommand{\thefootnote}{\arabic{footnote}}\footnotemark\footnotetext{The
stationarity condition of $f_t$ in Assumption 3 can be weakened to
Harris recurrence. See Bandi and Moloche (2004) for asymptotic
normality of local constant estimator under recurrence assumption. }
3--5, the state-domain estimator has the following asymptotic
normality
\begin{align*}
\sqrt{Nh}\mathrm{\
vech}\big(\widehat{\bSig}_{S,t}(x)-\bSig(x)-\frac{1}{2}h^2\mu_2\bSig''(x)\big)\toD
\mathcal{N}(0,2\nu_0p(x)^{-1} \Lambda(x)),
\end{align*}
where $\bSig''(x)$ is the matrix whose entries are the second
derivatives of the corresponding entries of $\bSig(x)$.
\end{theorem}

\medskip


Proposition \ref{P2} and Theorem \ref{T2} are both studied under the
assumption of high frequency data over a long time horizon, i.e.,
$\Delta\rightarrow0$ and $N\Delta\rightarrow\infty$. Various studies
under this assumption include Arfi (1998),
Gobet (2002), and Fan and Zhang (2003).

\bigskip

\section{DYNAMIC AGGREGATION OF TIME- AND STATE-DOMAIN ESTIMATORS}\label{sec3}

In this section, we show that the nonparametric estimators in the
time and state domains are asymptotically independent. This allows
us to combine these two estimators together to yield a more
efficient one.

\subsection{Asymptotic Normality}

The time- and state-domain estimators defined in the previous
sections are both driven by the factor process $f_t$. Intuitively,
with high probability, most of the data they use are far apart in
time. Since the Markov process $f_t$ is stationary and satisfies
some mixing condition (Assumption 3), it is reasonable to expect
that the time- and state-domain nonparametric estimators are also
asymptotically independent. The following theorem formally shows
this result. 

\begin{theorem}\label{T3}   Under the conditions of Theorems \ref{T1}
and \ref{T2}, conditioning on $f_t = x$, we have

(i) asymptotic independence:
\begin{eqnarray*}
\left( \begin{array}{c}  \sqrt{Nh}\ \mathrm{vech}\left(\widehat{\bSig}_{S,t}-\bSig(x)-\frac{1}{2}h^2\mu_2\bSig''(x)\right)\\
                          \sqrt{n}\ \mathrm{vech}\left(\widehat{\bSig}_{T,t}-\bSig(x)\right)
           \end{array} \right)\quad\quad\quad\quad\quad\quad\quad\quad \\
\quad\quad\quad\quad\quad\quad\quad\quad \stackrel{\mathcal
D}{\longrightarrow} {\mathcal N}\left(0,\left(\begin{array}{cc}
2\nu_0p(x)^{-1}\Lambda(x) & {\bf 0}
\\ {\bf
0}&\frac{\tau(1+e^\tau)}{(e^\tau-1)}\Lambda(x)
\end{array}\right)\right).
\end{eqnarray*}

\medskip
(ii) asymptotic normality of the aggregrated estimator
$\widehat{\bSig}_{A,t}(x)$ in (\ref{a1}):
\begin{align*}
\sqrt{Nh}\
\mathrm{vech}\Big(\widehat{\bSig}_{A,t}(x)-\bSig(x)-\frac{1}{2}h^2\omega_t(x)\mu_2\bSig''(x)\Big)
\toD \mathcal{N}(0,\Omega(x)),
\end{align*}
where $\Omega(x)=\big(2\omega_t^2(x)\nu_0p(x)^{-1}
+b(1-\omega_t(x))^2\frac{\tau(1+e^\tau)}{(e^\tau-1)}\big)\Lambda(x)$,
provided that $\lim Nh/n = b$ for some positive constant $b$ and
$h=O(N^{-1/5})$.
\end{theorem}

\medskip

Note that the nonparametric estimator in the time domain uses $n$
data points and the nonparametric estimator in the state domain
effectively uses the amount $O(Nh)$ of data.  The condition $\lim
Nh/n = b$ ensures that both estimators effectively use the same
amount (order) of data, which avoids the trivial case that either
the time domain or the state domain dominates the performance.

\subsection{Choice of the Dynamic Weight}

A natural question is how to choose the dynamic wight $\omega_t(x)$.
By Theorem 3(i) and (\ref{a6}), it is easy to see that for any
allocation vector $\ba$, the asymptotic optimal weight is
\begin{align}\label{d1}
\omega_t(x)=\frac{b\tau(1+e^\tau)p(x)}{2\nu_0(e^\tau-1)+b\tau(1+e^\tau)p(x)},
\end{align}
which is independent of $\ba$.  This choice\footnote{The optimal
choice of weight is proportional to the effective number of data
points used for the state-domain and time-domain smoothing. It
always outperforms the choice with $\omega_t=1$ (state-domain
estimator) or $\omega_t = 0$ (time-domain estimator).} also
optimizes the performance of the aggregated covariance estimator
$\widehat{\bSig}_{A,t}(x)$. Indeed, by Theorem 3(ii), the asymptotic
covariance matrix of $\widehat{\bSig}_{A,t}(x)$ is given by
$\Omega(x)$. It depends on the weight through the coefficient
$$\psi_t(x)\equiv2\omega_t^2(x)\nu_0p(x)^{-1}
+b(1-\omega_t(x))^2\frac{\tau(1+e^\tau)}{(e^\tau-1)},
$$
which is a quadratic function, and attains its minimum at
(\ref{d1}).

When $0<b<\infty$, the effective sample sizes in the time and state
domains are comparable. Hence, neither the time-domain nor the-state
domain estimator dominates. Therefore, by aggregating the time- and
state-domain estimators, we obtain an optimal reduction of
asymptotic variance.  The biases of the aggregated estimator are
indirectly controlled, when the optimal smoothing is conducted for
both time- and state-domain estimators so that their biases and
variances are already traded off before aggregation.

Note that at time $t$, the optimal weight $\omega_t(x)$ depends on
the current value of the factor process $f$ through the density
function $p(x)$. This is consistent with our common sense. When $f$
is low or high, $p(x)$ and consequently, the optimal weight are
approximately zero. In this case, the main contribution to the
aggregated estimator comes from the time-domain estimator. 
When $f$ is well in middle of its state space, say near its
unconditional mathematical expectation, the state-domain estimator
tends to dominate the aggregated estimator.

In practice, the density function $p(x)$ is unknown and should be
estimated. There are lots of existing methods to do this, such as
the kernel density estimator and the local time density estimator
(see A\"it-Sahalia, 1996; and Dalalyan and Kutoyants, 2003).


\bigskip

\section{NUMERICAL ANALYSIS}\label{sec4}

To evaluate the aggregated estimator, we compare it with the
time-domain estimator and the state-domain estimator. For the
time-domain estimation, we apply the exponential
smoothing\footnote{The choice comes from the recommendation of the
RiskMetrics of J.P. Morgan.  The parameter $\lambda$ can also be
chosen automatically by data by using the prediction error as in
Fan, Jiang, Zhang and Zhou (2003).  Since we compare the relative
performance between the time-domain estimator and the aggregated
estimator, we opt for this simple choice.  The results do not expect
to change much when a data-driven technique is used.} with
$\lambda=0.94$. For the state-domain estimation, we choose one yield
process as the ``factor," and then use it to estimate the volatility
matrix. The Epanechnikov kernel is used with the bandwidth $h$
chosen by generalized cross validation method (see Fan and Yao,
2003). To choose the optimal weight $\omega_t(x)$, we estimate the
density function $p(x)$ by the kernel density estimator (see
A\"it-Sahalia, 1996).

The following three measures are employed to assess the performance
of different methods for estimating the diffusion matrix. The first
two can only be used in simulation, and the last one can be used in
both simulation and real data analysis.

{\bf Measure 1.} The entropy loss is given by
\begin{align*} l_1(\bSig_t,\
\widehat{\bSig}_t)=\text{tr}(\bSig_t^{-1}\widehat{\bSig}_t)-\log|\bSig_t^{-1}\widehat{\bSig}_t|-\text{dim}(\bSig_t).
\end{align*}

{\bf Measure 2.} The quadratic loss is defined as
\begin{align*}
l_2(\bSig_t,\
\widehat{\bSig}_t)=\text{tr}\big(\widehat{\bSig}_t-\bSig_t\big)^2.
\end{align*}
%

{\bf Measure 3.} The prediction error (PE) is computed as
%
\begin{align}\label{e1}
\PE(\widehat{\bSig}_t)=\frac{1}{m}\sum_{i=T+1}^{T+m}\text{tr}\big(\bY_i\bY_i^T-\widehat{\bSig}_{t_i}\big)^2
\end{align}
for an out-sample of size $m$. The expected value can be decomposed
as
\begin{align*}
E[\PE(\widehat{\bSig}_t)]=\frac{1}{m}\sum_{i=T+1}^{T+m}E[\text{tr}\left(\bY_i\bY_i^T-\bSig_{t_i}\right)^2]+
\frac{1}{m}\sum_{i=T+1}^{T+m}E[\text{tr}\big(\bSig_{t_i}-\widehat{\bSig}_{t_i}\big)^2].
\end{align*}
Note that the second item reflects the effectiveness of the
estimated diffusion matrix, while the first term is the size of the
stochastic error, independent of the estimators. The first term is
usually an order of magnitude larger than the second term. Thus, a
small improvement in PE means a substantial improvement in estimated
volatility.  This will also be clearly demonstrated in our
simulation study (see Figure 4).

{\bf Measure 4.} Adaptive prediction error (APE).

As seen above, the dominant part of the PE is the stochastic error;
however, what we really care about is the estimation error. To
reduce the stochastic error in (\ref{e1}), we define the following
adaptive prediction error:
\begin{align}\label{e2}
\APE(\widehat{\bSig}_t)=\frac{1}{m}\sum_{i=T+1}^{T+m}\text{tr}\big(\frac{1}{2k+1}
\sum_{j=i-k}^{i+k}\bY_j\bY_j^T-\widehat{\bSig}_{t_i}\big)^2,
\end{align}
where $k$ is a nonnegative integer.  The basic idea is to average
out the stochastic errors first before computing square losses, but
this creates bias when $k$ is large. When $k=0$, the APE reduces to
the PE defined in (\ref{e1}).

\subsection{Simulation}

We use an essentially affine market price of risk specifications in
Duffee (2002) to simulate bond yields, and hence to obtain simulated
multivariate time series. 
Essentially affine model is the multivariate extension of the
square-root process. It has been proved useful in forecasting future
yields (see Duffee, 2002). Cheridito, Filipovi\'{c} and Kimmel
(2005) investigate the essentially affine model with one, two, and
three state variables, and give estimates of the parameters. We use
their one state variable model to conduct the simulations.


The one state variable affine term structure model assumes that the
instantaneous nominal interest rate $r_t$ is given by
$$
r_t=d_0+d_1s_t,
$$
where $d_0$ and $d_1$ are scalars, and $s_t$ is a scalar state
variable. The evolution of the state variable $s_t$ under the the
risk-neutral measure $Q$ is assumed to be
\begin{equation}\label{e2a}
ds_t=\big( a_1^Q+ b_{11}^Q s_t\big)dt+\sqrt{s_t} dW_t^Q.
\end{equation}
This is the well-known Cox-Ingersoll-Ross (CIR) model.

Let $P(t,\tau)$ be the time-$t$ price of a zero-coupon bond
maturing at $t+\tau$. Under the affine term structure and
the assumption of no arbitrage,
Duffie and Kan (1996)
show that the bond price admits the form
\begin{align}\label{e3}
P(t,\tau) = E_t^Q \exp(- \int_t^{t+\tau} r_u du )
=\exp[A(\tau)-B(\tau)s_t],
\end{align}
where $A(\tau)$ and $B(\tau)$ are both scalar functions satisfying
the following ordinary differential equations (ODEs)
\begin{equation}\label{e4}
\frac{dA(\tau)}{d\tau}=-a_1^QB(\tau)-d_0 \text{  and  }
\frac{dB(\tau)}{d\tau}=b_{11}^QB(\tau)-\frac{1}{2}B^2(\tau)+d_1.
\end{equation}
Thus, the bond's yield
\begin{equation}\label{e6}
y(s_t,\tau)=-\frac{1}{\tau}\log P(t,\tau)
=\frac{1}{\tau}[-A(\tau)+B(\tau)s_t]
\end{equation}
is affine in the state variable $s_t$.

We use the above model to simulate 5 zero-coupon bond yield
processes with maturities 1 month, 2 years, 4 years, 6 years, and 8
years. Since there is only one state variable $s_t$, the bond yields
of different maturities are perfectly linearly related, as shown in
(\ref{e6}), which is an unrealistic artifact of the model. To
attenuate this dilemma, Cherito \etal\ (2005) assume that only the
1-month yield process is observed without error, while other yields
are contaminated with i.i.d. multivariate Gaussian errors with mean
zero and unknown covariance matrix. They estimate the unknown
parameters from the yields of zero-coupon bonds extracted from the
US Treasury security prices from January 1972 to December 2002. The
estimated parameters are $a_1^Q=0.5$, $ b_{11}^Q=-0.0137$,
$d_0=0.0110$, and $ d_1=0.0074$. The standard deviations of the
Gaussian errors are estimated as $\sigma_1=0.0119,\
\sigma_2=0.0144,\ \sigma_3=0.0155$, and $ \sigma_4=0.0159$ for the
yields of 2-, 4-, 6-, and 8-year bonds, respectively. The associated
correlation coefficients are estimated as $\rho_{12}=0.9727$,
$\rho_{13}=0.9511$, $\rho_{14}=0.9371$, $\rho_{23}=0.9950$,
$\rho_{24}=0.9877$, and $\rho_{34}=0.9978$.

\begin{center}
Figure 2 here.
\end{center}

In the simulation, we set the the parameter values to be the above
estimated values from Cherito \etal\ (2005). We first generate
discrete samples of the state variable $s_t$ from diffusion process
(\ref{e2a}). Then we solve ODEs in (\ref{e4}) numerically. Figure 2
shows the solution to (\ref{e4}). After that, we obtain the ideal
yield processes by using (\ref{e6}) with maturities 1 month, 2
years, 4 years, 6 years, and 8 years. Finally, we add the i.i.d.
4-variate normal errors to the last 4 ideal yield
processes 
to obtain the observed bond processes with these
maturities\renewcommand{\thefootnote}{\arabic{footnote}}\footnotemark\footnotetext{
Here we add normal noise to make the model more realistic. Our
method performs even better without noise. Since the noise vectors
are i.i.d. across time and the standard deviations are small, adding
them to the original time series does not change the whole
structure. Hence, our theory can carry through under
contamination.}.

To generate the sample path of $s_t$, we use the transition density
property of the process. That is, given $s_t=x$, the variable
$2cs_{t+\Delta}$ has a noncentral chi-squared distribution with
degrees of freedom $4a_{11}^Q$ and noncentrality parameter
$2cxe^{b_{11}^Q\Delta}$, where
$c=\frac{2b_{11}^Q}{\exp(b_{11}^Q\Delta)-1}$. The initial value of
$s_0$ is generated from the invariant distribution of $s_t$, which
is gamma distribution with density
$p(y)=\frac{\omega^\nu}{\Gamma(\nu)}y^{\nu-1}e^{-\omega y}$, where
$\nu=2a_{11}^Q$ and $\omega=-2b_{11}^Q$.

We simulate 500 series of 1200 observations of weekly data with
$\Delta=1/52$ for the yields of five zero-coupon bonds with
maturities 1 month, 2 years, 4 years, 6 years, and 8 years,
respectively. For each simulated series, we set the last 150
observations as the out-sample data. 
For time $t$ out-sample data point, the time-domain estimator is
based on the past $n=104$ (two
years)\renewcommand{\thefootnote}{\arabic{footnote}}\footnotemark\footnotetext{With
$\lambda=0.94$, the last data point used in the time domain has an
extra weight $0.94^{104}\approx 0.0016$, which is very small. Hence,
we essentially include all the effective data points. }
observations, i.e., observations from $t-104$ to $t-1$;  and the
state-domain estimator is based on the 1050 data points right before
the current time, i.e.,
the data points from time $t-1050$ to $t-1$. 
The first yields process (1-month) is used as the factor for
state-domain estimation.

As pointed out in Section 1, the conditional covariance matrix of
the multivariate diffusion can be approximated by the diffusion
matrix times the sampling interval $\Delta$. Hence, we first obtain
estimates of the diffusion matrix, and then convert them into the
conditional covariance matrix estimates. The theoretical value of
the conditional variance of $s_t$ is given by Duffee (2002). Since
the bond yields are linear regression models of the state variable
(see (\ref{e6}) with Gaussian errors), the true (theoretical) value
of the conditional covariance matrix of the bond yields can be
easily obtained. By comparing the estimated conditional covariance
matrix to its theoretical value, the performance of our estimation
procedures is evaluated.

\begin{center}
Figure 3 here.
\end{center}

Figure 3 depicts the averages and standard deviations of the entropy
and quadratic losses of time-domain, state-domain, and aggregated
estimators. It shows unambiguously that the aggregated method always
has the smallest averages and standard deviations across 500
simulations for both the entropy loss and quadratic loss. Figures
4(a) and 4(b) summarize the distributions of the average losses over
150 out-samples forecasting across the 500 simulations. The results
are consistent with those in Figure 3. On the other hand, if the PE
in (\ref{e1}) with $m=150$ is used, the distributions look quite
different, which is demonstrated in Figure 4(c). It shows clearly
that even though there are huge efficiency improvements in
estimating the volatility matrix by using the aggregated method, the
improvements are masked by stochastic errors which are an order of
magnitude larger than the estimation errors. The average prediction
errors over 500 simulations are $1.850\times10^{-2}$,
$1.825\times10^{-2}$, and $1.846\times10^{-2}$ for the time-domain,
the aggregated, and the state-domain estimators, respectively. This
demonstrates that a small improvement in PE means a huge improvement
in the estimation of the volatility matrix. This effect is more
illuminatingly illustrated in Figure 4(d) where each point
represents a simulation.  The $x$-axis represents the ratios of the
averages of 150 quadratic losses for the time-domain estimator and
the state-domain estimator to those for the aggregated estimator,
whereas the $y$-axis is the ratios of the PEs for the time-domain
estimator and the state-domain estimator to those for the aggregated
estimator. The $x$-coordinates are mostly greater than 1, showing
the improved efficiency of the aggregated estimation. On the other
hand, the improved efficiency is masked by stochastic errors,
resulting in the $y$-coordinate spreading around the line $y=1$.

\begin{center}
Figure 4 here.
\end{center}

We have proved theoretically that nonparametric estimators based on
time-domain smoothing and state-domain smoothing are asymptotically
independent. To verify this, we compute their correlation
coefficients. Since both estimators are matrices, for a given
portfolio  allocation vector $\ba$, we compute the correlation of
the two estimators $\ba^T\widehat{\bSig}_{T,t}\ba$ and
$\ba^T\widehat{\bSig}_{S,t}\ba$ across 500 simulations at each given
time $t$ in the out-sample. Figure 5 presents the correlation
coefficients for $\ba=(0.2,\ 0.2,\ 0.2,\ 0.2,\ 0.2)^T$. Most of the
correlations are below 0.1, which strongly supports our theoretical
result. We also include the 95\% confidence intervals based on the
Fisher transformation in the same graph (the two dashed curves). A
large amount of these confidence intervals contain 0. The two
straight lines in the plot indicate the acceptance region for
testing the null hypothesis that the correlation coefficients are
zero at the significance level 5\%. Most of these null hypotheses
are accepted or nearly accepted. In fact, we conducted experiments
on the same simulations with larger sample sizes, and found that as
the sample size increases, the absolute values of the correlation
coefficients decrease to 0.

\begin{center}
Figure 5 here.
\end{center}

\subsection{Empirical Studies}

In this section, we apply the aggregated method to two sets of
financial data. Our aim is to examine whether our approach still
outperforms the time-domain and state-domain nonparametric
estimators in real applications.

\subsubsection{Treasury Bonds}

We consider the weekly returns of five treasury bonds with
maturities 3 months, 2 years, 5 years, 7 years, 10 years, and 30
years. We set the last 150 observations, which run from April 9,
1999 to February 15, 2002, as the out-sample data. For each
observation from the out-sample data, we use the past 104
observations (2 years) with $\lambda=0.94$ to obtain the time-domain
estimator, and the state-domain estimate is based on the past $900$
data points. The prediction error (Measure 3) and adaptive
prediction error (Measure 4) are used to assess the performance of
the three estimators: the time-domain estimator, the state-domain
estimator, and the aggregated estimator. The results are reported in
Table 1. From the table, we see that the aggregated estimator
outperforms significantly the other two estimators.

For comparison, the results from the simulated data are also
reported. Even through there is only a small improvement in PE for
simulated data, as evidenced in Section 4.1, there is a huge
improvement in the precision of estimating $\bSig_t$ in terms of
entropy loss (measure 1) and quadratic loss (measure 2).  Hence,
with the improvement of the PE in the bond price by the aggregrated
method, we would expect to have a huge improvement on the precision
of the estimation of covariance, which is of primary interest in
financial engineering.

\subsubsection{Exchange Rate}

We analyze the weekly exchange rates of five foreign currencies with
US dollars from September 6, 1985 to August 19, 2005. The five
foreign currencies are the Canadian Dollar, Australian Dollar,
Europe
Euro\renewcommand{\thefootnote}{\arabic{footnote}}\footnotemark\footnotetext{Europe
used several common currencies prior to the introduction of the
Euro. The European Currency Unit (ECU) was used from January 1, 1979
to January 1, 1999, when the Euro replaced the European Currency
Unit at par.}, UK British Pound, and Switzerland Franc. The length
of the time series is 1042. The exchange rates from December 6, 2002
to August 19, 2005, which are of length 142, are regarded as
out-sample data, and the estimation procedures are the same as
before, i.e., for each out-sample observation, the last 104 data
points with $\lambda = 0.94$ are set to construct the time-domain
estimator, the $900$ data points before the current time are used to
construct state-domain estimator, and then roll over. The results,
based on the PE and APE defined in Section 4, are also summarized in
Table 1. They demonstrate clearly that the aggregated estimator
outperforms the time-domain and state-domain estimators.

Using again the simulated data for calibration, as argued at the end
of Section 4.2.1, we would reasonably expect that the covariance
matrix estimated by the aggregated method outperforms significantly
both the matrices estimated by either the time- or state-domain
method alone.

\begin{center}
Table 1 here.
\end{center}

\bigskip

\section{DISCUSSIONS}

We have proposed an aggregated method to combine the information
from the time domain and state domain in multivariate volatility
estimation. To overcome the \emph{curse of dimensionality}, we
proposed a ``factor" modeling strategy. The performance comparisons
are studied both theoretically and empirically. We have shown that
the proposed aggregated method is more efficient than the estimators
based only on recent history or remote history. Our simulation and
empirical studies have also revealed that proper use of information
from both the time domain and the state domain makes volatility
matrix estimate more accurate. Our method exploits the continuity in
the time domain and stationarity in the state domain. It can also be
applied to situations where these two conditions hold approximately.

Our study has also revealed another potentially important
application of our method. It allows us to test the stationarity of
diffusion processes. When time-domain estimates differ substantially
from those of the state domain, it is an indication that the
processes is not stationary. Since the time-domain and state-domain
nonparametric estimators are asymptotically independent and normal,
formal tests can be formed. Further study on this topic
is beyond the scope of this paper.\\




\bigskip

\setcounter{equation}{0}
\renewcommand{\theequation}{A.\arabic{equation}}

\noindent{\Large\bf APPENDIX: PROOFS}

\bigskip

\noindent{\large\bf A.1\quad Proof of Proposition 1}

\medskip

\noindent In all the proofs below, we use $M$ to denote a generic
constant.

First, we show that the process $\{f_t\}$ is locally H\"{o}lder
continuous with order $q=(2q_0-1)/(4q_0)$ and coefficient $K_1$
satisfying $E[K_1^{4(q_0+\delta)}]<\infty$, i.e.
\begin{equation}\label{pf3}
|f_s-f_u|\leq K_1|s-u|^q,\quad s,u\in [t-\eta,t],
\end{equation}
where $\eta$ is a positive constant. Note that
\begin{align}\label{pf4}
\nonumber E|f_u-f_s|^{4(q_0+\delta)}& \leq ME\big|\int_s^u
a(f_v)dv\big|^{4(q_0+\delta)}+ ME\big|\int_s^u\sum_j
b_j(f_v)dW_v^j\big|^{4(q_0+\delta)}\\
&\equiv (I)+(II).
\end{align}
Then by Jensen's inequality and Assumption 2, we have
\begin{align}\label{pf5}
(I)\leq
M(u-s)^{4(q_0+\delta)-1}\int_s^uE|a(f_v)|^{4(q_0+\delta)}dv\leq
M(u-s)^{4(q_0+\delta)}.
\end{align}
On the other hand, applying martingale moment inequalities (see,
e.g. Karatzas and Shreve (1991), Section 3.3.D, p.163), Jensen's
inequality, and Assumption 2 gives
\begin{align}\label{pf6}
(II) \leq&  ME\big(\int_s^u\sum_j b_j^2(f_v)dv\big)^{2(q_0+\delta)}
\leq M(u-s)^{2(q_0+\delta)-1}\int_s^u\sum_j E|b_j(f_v)|^{4(q_0+\delta)}dv  \\
\nonumber\leq& M(u-s)^{2(q_0+\delta)}.
\end{align}
Combining (\ref{pf4}), (\ref{pf5}) and (\ref{pf6}) together leads to
$$
E|f_u-f_s|^{4(q_0+\delta)}\leq M(u-s)^{2(q_0+\delta)}.
$$
Thus by Theorem 2.1 of Revuz and Yor (1999, Page 26), we have
\begin{equation}\label{pf7}
E\big[\big(\sup_{s\neq
u}\{|f_s-f_u|/|s-u|^\alpha\}\big)^{4(q_0+\delta)}\big]<\infty
\end{equation}
for any $\alpha\in[0,\frac{2(q_0+\delta)-1}{4(q_0+\delta)})$. Let
$\alpha=\frac{2q_0-1}{4q_0}$ and $K_1=\sup_{s\neq
u}\{|f_s-f_u|/|s-u|^{\frac{2q_0-1}{4q_0}}\}$. Then
$E[K_1^{4(q_0+\delta)}<\infty]$, and inequality (\ref{pf3}) holds.

Second, by (\ref{pf1}) we have
$$\|\bsigma(f_s)-\bsigma(f_u)\|\le k_0|f_s-f_u|.$$
This together with (\ref{pf3}) shows that
\begin{eqnarray*}
\|\bsigma(f_s)-\bsigma(f_u)\|\leq k_0K_1|s-u|^q \equiv K|s-u|^q.
\end{eqnarray*}
Hence, $ E[K^{4(q_0+\delta)}]\leq M E[K_1^{4(q_0+\delta)}]< \infty$.
 $\blacksquare$

\bigskip
\bigskip

\noindent{\large\bf A.2\quad Proof of Theorem 1}

\medskip

\noindent{\bf Proof}. At time $s$, for fixed $k$, $\ell$, and $i$,
define $Z_{i,s}^{k,\ell}=(X_s^k-X_{t_i}^k)(X_s^\ell-X_{t_i}^\ell)$.
Applying Ito's formula to $Z_{i,s}^{k,\ell}$ results in
\begin{align*}
dZ_{i,s}^{k,\ell}=&(X_s^k-X_{t_i}^k)dX_s^\ell+(X_s^\ell-X_{t_i}^\ell)dX_s^k+
\sum_{j=1}^m\sigma_{kj}(f_s)\sigma_{\ell j}(f_s)ds\\
=&\left[(X_s^k-X_{t_i}^k)\mu_\ell(f_s)+(X_s^\ell-X_{t_i}^\ell)\mu_k(f_s)\right]ds\\
&+[\int_{t_i}^s\be_k^T\bmu(f_u)du\be_\ell^T\bsigma(f_s)+
\int_{t_i}^s\be_\ell^T\bmu(f_u)du\be_k^T\bsigma(f_s)]d\bW_s\\
&+[\int_{t_i}^s\be_k^T\bsigma(f_u)d\bW_u\be_\ell^T\bsigma(f_s)+\int_{t_i}^s
\be_\ell^T\bsigma(f_u)d\bW_u\be_k^T\bsigma(f_s)]d\bW_s\\
&+\sum_{j=1}^m\sigma_{kj}(f_s)\sigma_{\ell j}(f_s)ds.
\end{align*}
Hence, $Y_i^kY_i^\ell$ can be decomposed as
\begin{align*}
Y_i^kY_i^\ell=\Delta^{-1}Z_{i,t_{i+1}}^{k,\ell}\equiv
a_{i}^{k,\ell}+b_{i}^{k,\ell}+\bar{v}_{i}^{k,\ell},
\end{align*}
where
\begin{align*}
a_i^{k,\ell}=&\Delta^{-1}\int_{t_{i}}^{t_{i+1}}[(X_s^k-X_{t_i}^k)\mu_\ell(f_s)
+(X_s^\ell-X_{t_i}^\ell)\mu_k(f_s)]ds\\
&+\Delta^{-1}\int_{t_{i}}^{t_{i+1}}\int_{t_i}^s[\be_k^T\bmu(f_u)du\be_\ell^T\bsigma(f_s)+
\be_\ell^T\bmu(f_u)du\be_k^T\bsigma(f_s)]d\bW_s,
\end{align*}
\begin{align*}
b_i^{k,\ell}=\Delta^{-1}\int_{t_i}^{t_{i+1}}\int_{t_i}^s[\be_k^T\bsigma(f_u)d\bW_u\be_\ell^T\bsigma(f_s)+
\be_\ell^T\bsigma(f_u)d\bW_u\be_k^T\bsigma(f_s)]d\bW_s
\end{align*}
and
\begin{align*}
c_i^{k,\ell}=\Delta^{-1}\int_{t_{i}}^{t_{i+1}}\sum_{j=1}^m\sigma_{kj}(f_s)\sigma_{\ell
j}(f_s)ds.
\end{align*}
Correspondingly, (\ref{b4}) has the following decomposition
\begin{align}
\nonumber\widehat{U}_{T,t}&=\frac{1-\lambda}{1-\lambda^n}
\sum_{i=1}^n\lambda^{i-1}\sum_{\ell\leq k}c_{k\ell}a_{N-i}^{k,\ell}+
\frac{1-\lambda}{1-\lambda^n}
\sum_{i=1}^n\lambda^{i-1}\sum_{\ell\leq
k}c_{k\ell}b_{N-i}^{k,\ell}\\
\nonumber&\quad\quad+ \frac{1-\lambda}{1-\lambda^n}
\sum_{i=1}^n\lambda^{i-1}\sum_{\ell\leq k}^kc_{k\ell}\bar{v}_{N-i}^{k,\ell}\\
&\equiv A_{n,\Delta}+B_{n,\Delta}+V_{n,\Delta}.\label{005}
\end{align}
Therefore, Slutsky's lemma, together with Lemmas \ref{L1}--\ref{L3}
below, leads to the conclusions of Theorem \ref{T1} immediately.
 $\blacksquare$

\begin{lemma} \label{L1}  Under Assumption 1, as
$n\rightarrow\infty,\ n\Delta\rightarrow0$, and
$n(1-\lambda)\rightarrow \tau$, we have
\begin{equation}\label{021}
EA_{n,\Delta}^2=O(\Delta),
\end{equation}
where $A_{n,\Delta}=\frac{1-\lambda}{1-\lambda^n}
\sum_{i=1}^n\lambda^{i-1}\sum_{\ell\leq
k}c_{k\ell}a_{N-i}^{k,\ell}$, as defined in (\ref{005}).
\end{lemma}


\noindent{\bf\em Proof}. First, note that
\begin{align}
 E(a_{i}^{k,\ell})^2\leq&\label{pf20}
2E\big(\Delta^{-1}\int_{t_{i}}^{t_{i+1}}[(X_s^k-X_{t_i}^k)\mu_\ell(f_s)
+(X_s^\ell-X_{t_i}^\ell)\mu_k(f_s)]ds\big)^2\\
\nonumber&+2E\big(\Delta^{-1}\int_{t_i}^{t_{i+1}}\int_{t_i}^s[\be_k^T\bmu(f_u)du\be_\ell^T\bsigma(f_s)+
\be_\ell^T\bmu(f_u)du\be_k^T\bsigma(f_s)]d\bW_s\big)^2\\
\nonumber\equiv& I_1(\Delta)+I_2(\Delta).
\end{align}

Applying Jensen's inequality and H\"{o}lder's inequality (Propostion
1), we obtain
\begin{align}
 I_1(\Delta)\leq&
M\Delta^{-1}\int_{t_i}^{t_{i+1}}E\left[(X_s^k-X_{t_i}^k)\mu_\ell(f_s)
+(X_s^\ell-X_{t_i}^\ell)\mu_k(f_s)\right]^2ds\label{pf8}\\
\nonumber\leq&M\Delta^{-1}\int_{t_i}^{t_{i+1}}\Big\{\big(E(X_s^k-X_{t_i}^k)^4E[\mu_\ell(f_s)]^4\big)^{1/2}
+\big(E(X_s^\ell-X_{t_i}^\ell)^4E[\mu_k(f_s)]^4\big)^{1/2}\Big\}ds.
\end{align}
Since an application of Jensen's inequality, martingale moments
inequalities and Assumption 2 results in
\begin{align*}
E(X_s^\ell-X_{t_i}^\ell)^4&\leq
M\big(E\big[\int_{t_i}^s\mu_\ell(f_u)du\big]^4+
\sum_{j=1}^mE\big[\int_{t_i}^s\sigma_{\ell j}(f_u)dW_u^j \big]^4\big)\\
&\leq M\big((s-t_i)^3\int_{t_i}^sE[\mu_\ell(f_u)]^4du+
\sum_{j=1}^mM(s-t_i)\int_{t_i}^sE[\sigma_{\ell j}(f_u)]^4du\big)\\
&\leq M(s-t_i)^2,
\end{align*}
we see that (\ref{pf8}) can be bounded as
\begin{align}\label{pf9}
I_1(\Delta)\leq M\Delta.
\end{align}

We now consider the second term $I_2(\Delta)$ in (\ref{pf20}). By
stochastic calculus and Jensen's inequality, we have
\begin{align*}
I_2(\Delta)&=2\int_{t_i}^{t_{i+1}}\sum_{j=1}^mE\Big(\Delta^{-1}\int_{t_i}^s[\mu_k(f_u)\sigma_{\ell
j}(f_s)+ \mu_\ell(f_u)\sigma_{kj}(f_s)]du\Big)^2ds\\
&\leq
M\Delta^{-1}\int_{t_i}^{t_{i+1}}\sum_{j=1}^m\int_{t_i}^sE[\mu_k(f_u)\sigma_{\ell
j}(f_s)+
\mu_\ell(f_u)\sigma_{kj}(f_s)]^2duds\\
&=O(\Delta).
\end{align*}
This together with (\ref{pf9}) leads to
$E(a_{i}^{k,\ell})^2=O(\Delta) $. Therefore, by the Cauchy-Schwarz
inequality and the assumption that $\lim_n (1 - \lambda)$ exists,
\begin{align*}
EA_{n,\Delta}^2\leq Mn\left(\frac{1-\lambda}{1-\lambda^n}\right)^2
\sum_{i=1}^n\lambda^{2(i-1)}\sum_{\ell\leq
k}c_{k\ell}^2E(a_{N-i}^{k,\ell})^2= O(\Delta),
\end{align*}
which concludes the proof.  $\blacksquare$

\medskip

\begin{lemma} \label{L2} Under Assumptions 1 and 2, as
$n\rightarrow\infty,\  n\Delta^{q}\rightarrow0$ and
$n(1-\lambda)\rightarrow \tau$, we have
\begin{align*}
\sqrt{n}B_{n,\Delta}\toD Z_{\bc},
\end{align*}
where $B_{n,\Delta}$ is defined in (\ref{005}) and the random
variable $Z_{\bc}$ is defined in Theorem 1.
\end{lemma}

\noindent{\bf\em Proof}. We will decompose $B_{n,\Delta}$ into two
parts and prove that the first part is asymptotically negligible and
the second part has some asymptotic distribution.

Note that $b_{i}^{k,\ell}$ can be decomposed as
\begin{align}\label{004}
b_{i}^{k,\ell}=\mathcal{B}_i^{k,\ell}+ \mathcal{C}_i^{k,\ell},
\end{align}
where
\begin{align*}
\mathcal{B}_i^{k,\ell}=\Delta^{-1}\sum_{j,p}\left(\sigma_{kj}(f_{t_0})\sigma_{\ell
p}(f_{t_0})+\sigma_{kp}(f_{t_0})\sigma_{\ell
j}(f_{t_0})\right)\int_{t_i}^{t_{i+1}}(W_{s}^j-W_{t_i}^j)dW_s^p
\end{align*}
and
\begin{align*}
\mathcal{C}_i^{k,\ell}=\Delta^{-1}\int_{t_i}^{t_{i+1}}\int_{t_i}^s
[\be_k^T(\bsigma(f_u)-\bsigma(f_{t_0}))d\bW_u\be_\ell^T\bsigma(f_s)+
\be_k^T\bsigma(f_u)d\bW_u\be_\ell^T(\bsigma(f_s)-\bsigma(f_{t_0}))]d\bW_s,
\end{align*}
where $\be_k$ is the unit vector with $k$th entry 1 and all other
entries 0. Correspondingly, $B_{n,\Delta}$ is decomposed as
$$
B_{n,\Delta}=\frac{1-\lambda}{1-\lambda^n}\sum_{k\leq\ell}c_{k\ell}\sum\lambda^{i-1}\mathcal{B}_{N-i}^{k,\ell}
+\frac{1-\lambda}{1-\lambda^n}\sum_{k\leq\ell}c_{k\ell}\sum\lambda^{i-1}\mathcal{C}_{N-i}^{k,\ell}\equiv
\mathcal{B}+\mathcal{C}.
$$

First, we show that $\sqrt{n}\mathcal{C}$ is asymptotically
negligible. To this end, note that by stochastic calculus and the
triangular inequality, we have
\begin{align}
\nonumber E(\mathcal{C}_i^{k,\ell})^2\leq&
\Delta^{-2}\int_{t_i}^{t_{i+1}}\sum_{j=1}^m E\Big(\int_{t_i}^s
\be_k^T(\bsigma(f_u)-\bsigma(f_{t_0}))d\bW_u\sigma_{\ell
j}(f_s)\Big)^2ds\\
\nonumber&+\Delta^{-2}\int_{t_i}^{t_{i+1}}\sum_{j=1}^mE\Big(\int_{t_i}^s
\be_k^T\bsigma(f_u)d\bW_u(\sigma_{\ell j}(f_s)-\sigma_{\ell
j}(f_{t_0}))\Big)^2ds\\
\nonumber\equiv&
\Delta^{-2}\int_{t_i}^{t_{i+1}}\sum_{j=1}^mI_1^{(j)}(\Delta)ds
+\Delta^{-2}\int_{t_i}^{t_{i+1}}\sum_{j=1}^mI_1^{(j)}(\Delta)ds.
\end{align}
Applying H\"{o}lder's inequality yields
\begin{align}\label{029}
I_1^{(j)}(\Delta) &\leq \Big(E\big(\int_{t_i}^s
\be_k^T(\bsigma(f_u)-\bsigma(f_{t_0}))d\bW_u\big)^4E(\sigma_{\ell
j}(f_s))^4\Big)^{1/2},
\end{align}
and then by martingale moment inequalities and (\ref{b5}) we obtain
\begin{align*}
E\Big(\int_{t_i}^s
\be_k^T(\bsigma(f_u)-\bsigma(f_{t_0}))d\bW_u\Big)^4\leq& O(1)
E\Big(\int_{t_i}^s\sum_{j=1}^m(\sigma_{kj}(f_u)-\sigma_{kj}(f_{t_0}))^2du
\Big)^2\\
\leq& O\left((n\Delta+\Delta)^{4q}\Delta^2\right).
\end{align*}
Hence, we can bound (\ref{029}) as
\begin{align}\label{031}
I_1^{(j)}(\Delta)\leq O\left((n\Delta)^{2q}\Delta\right).
\end{align}
Next we consider $I_2^{(j)}(\Delta)$. Similarly, by H\"{o}lder's
inequalities, martingale moments inequalities, and (\ref{b5}) we
have
\begin{align*}
I_2^{(j)}(\Delta)&\leq\Big(E\big(\int_{t_i}^s
\be_k^T\bsigma(f_u)d\bW_u\big)^4E(\sigma_{\ell j}(f_s)-\sigma_{\ell
j}(f_{t_0}))^4\Big)^{1/2}\\
&\leq O(1)\Big( E[\int_{t_i}^s\sum_{j=1}^m\sigma_{kj}^2(f_u)du
]^2 (n\Delta+\Delta)^{4q}EK^4\Big)^{1/2}\\
&\leq O\left((n\Delta)^{2q}\Delta\right).
\end{align*}
This together with (\ref{031}) implies that
\begin{align*}
E(\mathcal{C}_i^{k,\ell})^2= O\left((n\Delta)^{2q}\right).
\end{align*}
Hence, it follows that
 \begin{equation}\label{034}
 E\left(\sqrt{n}\mathcal{C}\right)^2= O((n\Delta)^{2q}),
 \end{equation}
which means that $\sqrt{n}\mathcal{C}$ is asymptotically negligible.

Next, we consider the term $\sqrt{n}\mathcal{B}$.  We first define
the augmented filtration $\mathcal{F}_t$. Let $(\Omega, \mathcal{F},
P)$ be the probability space in which the Brownian motion $\{\bW_t,
0\leq t<\infty \}$ is defined, and $\bX_0$ is the initial value of
model (\ref{a2}) and independent of $\mathcal{F}_{\infty}$. Define
the left-continuous filtration
$\mathcal{G}_t=\sigma(\bX_0)\vee\{\mathcal{F}_t^W,0\leq t< \infty\}$
as well as the collection of null sets $\mathcal{N}=\{N\in\Omega;
\exists G\in \mathcal{G}_{\infty} \text{ with } N\subseteq G \text{
and } P(G)=0\}$. Then the augmented filtration is defined as
$\mathcal{F}_t=\sigma(\mathcal{G}_t\cup \mathcal{N}),\ 0\leq
t<\infty$;
$\mathcal{F}_{\infty}=\sigma(\bigcup_{t\geq0}\mathcal{F}_t)$. First
note that by stochastic calculus we have
$E[\mathcal{B}_i^{k,\ell}|\mathcal{F}_0]=0$ and for $i\neq j$,
$\mathcal{B}_i^{k,\ell}$ and $\mathcal{B}_j^{k,\ell}$ are
independent. Therefore, we only need to verify the conditions of the
central limit theorem for the martingale difference array (see, e.g.
Hall and Heyde (1980), Corollary 3.1, P.58); namely, we need to
check
\begin{align}\label{pf10}
\sum_{i=1}^nE\Big(\frac{\sqrt{n}(1-\lambda)}{1-\lambda^n}\lambda^{i-1}\sum_{\ell\leq
k}c_{k\ell}\mathcal{B}_{i}^{k,\ell}|\mathcal{F}_{t_{i}}\Big)^2
\toP\frac{\tau(1+e^\tau)}{e^\tau-1}\bc^TP_D^T(\bSig(f_t)\otimes\bSig(f_t))P_D\bc
\end{align}
and
\begin{align}\label{pf7a}
\sum_{i=1}^nE\Big[\big(\sqrt{n}\frac{1-\lambda}{1-\lambda^n}\lambda^{i-1}\sum_{\ell\leq
k}c_{k\ell}\mathcal{B}_{i}^{k,\ell}\big)^4\Big|\mathcal{F}_{t_i}\Big]\toP0.
\end{align}
Expression(\ref{pf10}) gives the asymptotic conditional variance of
$\sqrt{n}\mathcal{B}$ and (\ref{pf7a}) implies the conditional
Lindeberg condition. These two conditions lead to
\begin{align}\label{035}
\sqrt{n}\mathcal{B}\toD Z_{\bc},
\end{align}
where the random variable $Z_{\bc}$ is defined as in Theorem
\ref{T1}.

We first prove (\ref{pf10}). From stochastic calculus we know that
$E[\mathcal{B}_i^{k,\ell}|\mathcal{F}_{t_i}]=0$ and for $i\neq j$,
$\mathcal{B}_i^{k,\ell}$ and $\mathcal{B}_j^{k,\ell}$ are
independent. Moreover, by (\ref{b5}) we have
\begin{align*}
E[\mathcal{B}_{i}^{k_1,\ell_1}\mathcal{B}_{i}^{k_2,\ell_2}|\mathcal{F}_{t_{i}}]=&\Delta^{-2}\sum_{j,g}H_{j,g}^{k_1,\ell_1}(f_{t_0})
H_{j,g}^{k_2,\ell_2}(f_{t_0})\int_{t_i}^{t_{i+1}}E(W_{s}^j-W_{t_i}^j)^2ds\\
=&\frac{1}{2}\sum_{j,g}H_{j,g}^{k_1,\ell_1}(f_{t_0})
H_{j,g}^{k_2,\ell_2}(f_{t_0})\\
=&\frac{1}{2}\sum_{j,g}H_{j,g}^{k_1,\ell_1}(f_{t})
H_{j,g}^{k_2,\ell_2}(f_{t})+o_g((n\Delta+\Delta)^q),
\end{align*}
where $H_{j,g}^{k,\ell}(x)=\sigma_{kj}(x)\sigma_{\ell
g}(x)+\sigma_{kg}(x)\sigma_{\ell j}(x)$. It follows that
\begin{align*}
\var(\sum_{\ell\leq k}c_{\ell
k}\mathcal{B}_i^{\ell,k}|\mathcal{F}_{t_{i}})&=\bc^TP_D(2\bSig(f_{t_0})\otimes\bSig(f_{t_0}))P_D^T\bc\\
&\toP \bc^TP_D(2\bSig(f_{t})\otimes\bSig(f_{t}))P_D^T\bc.
\end{align*}
Therefore, we get the following result for the conditional variance
of the left hand side of (\ref{pf10}):
\begin{align*}
\sum_{i=1}^nE\Big(\frac{\sqrt{n}(1-\lambda)}{1-\lambda^n}\lambda^{i-1}\sum_{\ell\leq
k}c_{k\ell}\mathcal{B}_{i}^{k,\ell}|\mathcal{F}_{t_{i}}\Big)^2
&=\frac{n(1-\lambda)(1+\lambda^n)}{(1+\lambda)(1-\lambda^n)}\var(\sum_{\ell\leq
k}c_{\ell
k}\mathcal{B}_i^{\ell,k}|\mathcal{F}_{t_{i}})\\
&\toP\frac{\tau(1+e^\tau)}{e^\tau-1}\bc^TP_D^T(\bSig(f_t)\otimes\bSig(f_t))P_D\bc,
\end{align*}
where $\tau=\lim_{n\rightarrow\infty} n(1-\lambda)$. This verifies
(\ref{pf10}).

Then we show (\ref{pf7}). Straightforward calculations yield
\begin{align*}
&E\Big[\big(\sum_{\ell\leq
k}c_{k\ell}\mathcal{B}_{i}^{k,\ell}\big)^4\big|\mathcal{F}_{t_{i}}\Big]=
O(1)\sum_{\ell\leq k}c_{k\ell}^4E[(\mathcal{B}_{i}^{k,\ell})^4|\mathcal{F}_{t_{i}}]\\
=&O(1)\sum_{\ell\leq
k}c_{k\ell}^4\Delta^{-4}\sum_{j,g}(H_{j,g}^{k,\ell}(f_{t_0}))^4
E\Big[\big(\int_{t_i}^{t_{i+1}}(W_{s}^j-W_{t_i}^j)dW_s^g\big)^4\Big|\mathcal{F}_{t_{i-1}}\Big]\\
=&O(1)\sum_{\ell\leq
k}c_{k\ell}^4\sum_{j,g}(H_{j,g}^{k,\ell}(f_{t_0}))^4.
\end{align*}
This together with Assumption 2 and H\"{o}lder's inequality leads to
\begin{align*}
\sum_{i=1}^nE\Big[\big(\sqrt{n}\frac{1-\lambda}{1-\lambda^n}\lambda^{i-1}\sum_{\ell\leq
k}c_{k\ell}\mathcal{B}_{i}^{k,\ell}\big)^4\big|\mathcal{F}_{t_i}\Big]=
O(n^{-1})\sum_{\ell\leq
k}c_{k\ell}^4\sum_{j,g}(H_{j,g}^{k,\ell}(f_{t_0}))^4\toP0,
\end{align*}
which proves (\ref{pf7}). (\ref{035}) holds in consequence.
Combining (\ref{034}) and (\ref{035}) and applying Slutsky's lemma,
we obtain the conclusion in lemma \ref{L2}.
 $\blacksquare$

\medskip

\begin{lemma} \label{L3} Under Assumptions 1 and 2, as
$n\rightarrow\infty$ and $n\Delta^q\rightarrow0$, the following
result holds for $C_{n,\Delta}$ defined in (\ref{005})
\begin{equation}\label{022}
E\left|C_{n,\Delta}- \bc^T\vech(\bSig(f_t))\right|=
O\left((n\Delta)^q\right).
\end{equation}
\end{lemma}

\noindent{\bf\em Proof}. Note that
\begin{align*}
 E|C_{n,\Delta}- \sum_{\ell\leq k}^kc_{k\ell}v_{k\ell,t}|&=
\frac{1-\lambda}{1-\lambda^n}
E\big|\sum_{i=1}^n\lambda^{i-1}\sum_{\ell\leq
k}^kc_{k\ell}\left(\bar{v}_{N-i}^{k,\ell}- v_{k\ell,t}\right)\big|\\
&\leq \frac{1-\lambda}{1-\lambda^n}
\sum_{i=1}^n\lambda^{i-1}\sum_{\ell\leq
k}^kc_{k\ell}E|\bar{v}_{N-i}^{k,\ell}- v_{k\ell,t}|.
\end{align*}
Thus we only need to consider the asymptotic property of
$E|\bar{v}_{i}^{k,\ell}- v_{k\ell,t}|$. By the Cauchy-Schwarz
inequality and H\"{o}lder's inequality, we have
\begin{align*}
E\big|\bar{v}_{i}^{k,\ell}- v_{k\ell,t}\big|&\leq
\Delta^{-1}\sum_{j=1}^m\int_{t_{i}}^{t_{i+1}}\big\{E\big|\sigma_{kj}(f_t)\bigl(\sigma_{\ell
j}(f_t)-
\sigma_{\ell j}(f_s)\bigr)\big|\\
&\quad+E\big|
\bigl(\sigma_{kj}(f_t)-\sigma_{kj}(f_s)\bigr)\sigma_{\ell j}(f_s)\big|\big\}ds\\
&\leq
\Delta^{-1}\sum_{j=1}^m\int_{t_{i}}^{t_{i+1}}\big\{\big[E\sigma_{kj}^2(f_t)E\bigl(\sigma_{\ell
j}(f_t)-
\sigma_{\ell j}(f_s)\bigr)^2\big]^{1/2}\\
&\quad+\big[E
\bigl(\sigma_{kj}(f_t)-\sigma_{kj}(f_s)\bigr)^2E\sigma_{\ell
j}^2(f_s)\big]^{1/2}\big\}ds
\end{align*}
Therefore by (\ref{b5}) and Assumption 2,
\begin{align*}
E\big|\bar{v}_{i}^{k,\ell}-
v_{k\ell,t}\big|=O\big((n\Delta+\Delta)^q\big)=
O\big((n\Delta)^q\big).
\end{align*}
This proves (\ref{022}).  $\blacksquare$

\bigskip
\bigskip

\noindent{\large\bf A.3\quad Proof of Proposition 2}

\medskip

\begin{lemma} \label{L8}  Since $f_t$ is a stationary real ergodic
process, we have
\begin{equation*}
\frac{L_f(T,x)}{\sum b_j^2(x)T}\toas p(x),
\end{equation*}
where $p(x)$ is the time-invariant density function of the process
$f_t$ at $x$.
\end{lemma}

\noindent{\bf\em Proof}. See Bandi and Phillips (2003) and Bosq
(1998, Theorem 6.3, P150). $\blacksquare$

\bigskip


\begin{lemma} \label{L9}  Suppose $\Delta\rightarrow0$,
$N\Delta\rightarrow\infty$, and
$\frac{1}{h}\sqrt{\Delta\log\Delta^{-1}}=o(1)$. Under Assumptions
3--5, we have for $\ell=0,1,2,3$
$$
W_\ell(x)=\frac{1}{\Delta}\int_{t_0}^{t_{N-1}}(f_s -x)^\ell
K_h(f_s-x)ds+Nh^{\ell-1}O_{a.s.}\Bigl(\sqrt{\Delta\log\Delta^{-1}}\Bigr).$$
\end{lemma}

\noindent{\bf\em Proof}. First, note that for any nonnegative
integer $\ell\leq4$, we have
\begin{align}\label{pf11}
&\Big|W_\ell(x)-\frac{1}{\Delta}\sum_{k=0}^{N-1}\int_{t_k}^{t_{k+1}}(f_s
-x)^\ell K\Bigl(\frac{f_s-x}{h}\Bigr)ds\Big| \\
\nonumber\leq&\frac{1}{h\Delta}\sum_{k=0}^{N-1}\int_{t_k}^{t_{k+1}}\Bigl|(f_{t_k}
-x)^\ell K\Bigl(\frac{f_{t_k}-x}{h}\Bigr)-(f_s -x)^\ell
K\Bigl(\frac{f_s-x}{h}\Bigr)\Bigr|ds\\
\nonumber\leq& I_1+I_2
\end{align}
with
\begin{align}
I_1=\frac{1}{h\Delta}\sum_{k=0}^{N-1}\int_{t_k}^{t_{k+1}}\left|
K'\Bigl(\frac{\widehat{r}_{ks}-x}{h}\Bigr)\right|\left|\frac{f_s-f_{t_k}}{h}\right||f_{t_k}
-x|^\ell ds \label{pf12}
\end{align}
and
\begin{align}
I_2=\frac{1}{h\Delta}\sum_{k=0}^{N-1}\int_{t_k}^{t_{k+1}}\left|(\overline{r}_{ks}
-x)^{\ell-1}(f_s
-f_{t_k})\right|K\Bigl(\frac{f_s-x}{h}\Bigr)ds,\label{pf13}
\end{align}
where $\widehat{r}_{ks}$ and $\overline{r}_{ks}$ are both values on
the line segment connecting $f_{t_k}$ to $f_s$. Now define
$$
\kappa_{N,\Delta}=\max_{i\leq N-1}\sup_{t_{i-1}\leq s\leq
t_i}|f_s-f_{t_{i-1}}|.
$$
Then, by Levy's modulus of continuity of diffusions (see, e.g. Revuz
and Yor (1998, Ch. V, Exercise 1.20)),
\begin{align}\label{pf14}
P\Bigl(\limsup_{\Delta\rightarrow0}\frac{\kappa_{N,\Delta}}{\sqrt{\Delta\log\Delta^{-1}}}=\alpha
\Bigr)=1,
\end{align}
where $\alpha$ is a suitable constant. In turn, (\ref{pf14}) implies
that
\begin{align*}
\kappa_{N,\Delta}=O_{a.s.}\Bigl(\sqrt{\Delta\log\Delta^{-1}}\Bigr).
\end{align*}
This together with the assumption that
$\frac{1}{h}\sqrt{\Delta\log\Delta^{-1}}=o(1)$ leads to
\begin{align*}
\frac{\kappa_{N,\Delta}}{h}=o_{a.s.}(1) \text{ as }
N\Delta\rightarrow\infty.
\end{align*}
In view of (\ref{pf12}) and (\ref{pf13}), we have
\begin{align*}
K'\bigl(\frac{\widehat{r}_{ks}-x}{h}\bigr)=K'\bigl(\frac{f_s-x}{h}+o_{a.s.}(1)\bigr)
\end{align*}
and
\begin{align*}
\overline{r}_{ks} -x=h\big(\frac{f_s-x}{h}+o_{a.s.}(1)\big),
\end{align*}
uniformly over $k=0,\cdots,N-1$. Hence, by Lemma \ref{L8} and Revuz
and Yor (1999), Exercise 1.15 and Corollary 1.6 of Chapter 6, we
obtain that (\ref{pf12}) can be bounded as
\begin{align*}
I_1\leq&\frac{\kappa_{N,\Delta}}{h}\frac{h^{\ell-1}}{\Delta}\sum_{k=0}^{N-1}\int_{t_k}^{t_{k+1}}\big|
K'\bigl(\frac{f_s-x}{h}+o_{a.s.}(1)\bigr)\big|\big|\frac{f_s
-x}{h}+o_{a.s.}(1)\big|^\ell ds \\
=&N\Delta
h^{\ell-1}\frac{\kappa_{N,\Delta}}{h}\int_{-\infty}^{\infty}\big|
K'\bigl(\frac{y-x}{h}+o_{a.s.}(1)\bigr)\big|\big|\frac{y
-x}{h}+o_{a.s.}(1)\big|^\ell\frac{L_r(t_{N-1},y)}{N\Delta\sum
b_j^2(y)} dy\\
=&N h^{\ell}\frac{\kappa_{N,\Delta}}{h}\int_{-\infty}^{\infty}\big|
K'\bigl(u+o_{a.s.}(1)\bigr)\big||u+o_{a.s.}(1)|^\ell(p(uh+x)+o_{a.s.}(1))
du.
\end{align*}
This together with (\ref{pf2}) yields
$$
I_1\leq N
h^{\ell}O_{a.s.}\bigl(\frac{1}{h}\sqrt{\Delta\log\Delta^{-1}}\bigr).
$$
Similarly, we can show that (\ref{pf13}) is also bounded by
$Nh^{\ell}O_{a.s.}\bigl(\frac{1}{h}\sqrt{\Delta\log\Delta^{-1}}\bigr)$.
This proves the stated results.  $\blacksquare$

\bigskip

\noindent{\bf Proof of Proposition 2}

\medskip

\noindent Since $x^{2\ell}K(x)$ is a positive function, by Exercise
1.15 and Corollary 1.6 of Chapter 6 of Revuz and Yor (1999),  and
Lemma \ref{L8} above we have for $\ell=0,1$,
\begin{align*}
&\frac{1}{N\Delta}\int_{t_0}^{t_{N-1}}\big(\frac{f_s-x}{h}\big)^{2\ell}
K\bigl(\frac{f_s-x}{h}\bigr)ds
\\
&=\int\big(\frac{y-x}{h}\big)^{2\ell}
K\bigl(\frac{y-x}{h}\bigr)\frac{L_r(t_{N-1},y)}{N\Delta\sum b_j^2(y)} dy\\
&=h\int u^{2\ell} K(u)(p(uh+x)+o_{a.s.}(1))du\\
&=h\bigl(p(x)\mu_{2\ell}+o_{a.s.}(1)\bigr),
\end{align*}
where we have used $\mu_4=\int x^4K(x)dx<\infty$. This together with
Lemma \ref{L9} leads to
\begin{align}\label{pf15}
\frac{1}{N}W_{2\ell}(x)=&\frac{1}{N\Delta}\int_{t_0}^{t_{N-1}} (f_s
-x)^{2\ell}
K_h(f_s-x)ds+o_{a.s.}(1)\\
\nonumber=&h^{2\ell}(p(x)\mu_{2\ell}+o_{a.s.}(1)).
\end{align}

Let $s(dx)=\exp\Bigl\{\int_\alpha^x\frac{2a(y)}{\sum
b_j^2(y)}dy\Bigr\}\frac{2dx}{\sum b_j^2(x)}$ be the speed measure of
$f_t$. By the Quotient theorem (Revuz and Yor (1999), Theorem 3.12,
Chapter 10, p.427),
\begin{align*}
\frac{\frac{1}{N\Delta}\int_{t_0}^{t_{N-1}}\left(\frac{f_s
-x}{h}\right)^{2\ell+1}
K_h(f_s-x)ds}{\frac{1}{N\Delta}\int_{t_0}^{t_{N-1}}
K_h(f_s-x)ds}&=\frac{\int\left(\frac{y-x}{h}\right)^{2\ell+1}K_h(y-x)s(dy)}{
\int
K_h(y-x)s(dy)}+o_{a.s.}(1)\\
&=\frac{\mu_{2\ell+1}}{\mu_0}+o_{a.s.}(1)
\end{align*}
as $N\Delta\rightarrow\infty$. In turn, this implies that
 \begin{align}\label{pf16}
\frac{W_{2\ell+1}(x)/h^{2\ell+1}}{W_0(x)}&=\frac{\frac{1}{\Delta}\int_{t_0}^{t_{N-1}}\left(\frac{f_s
-x}{h}\right)^{2\ell+1}
K_h(f_s-x)ds+NO_{a.s.}\bigl(\frac{\sqrt{\Delta\log\Delta^{-1}}}{h}\bigr)}{\frac{1}{\Delta}\int_{t_0}^{t_{N-1}}
K_h(f_s-x)ds+NO_{a.s.}\bigl(\frac{\sqrt{\Delta\log\Delta^{-1}}}{h}\bigr)}\\
\nonumber&= \frac{\mu_{2\ell+1}}{\mu_0}+o_{a.s.}(1).
\end{align}
Combining (\ref{pf15}) and (\ref{pf16}), we obtain
\begin{align*}
W_{2\ell+1}(x)=Nh^{2\ell+1}(p(x)\mu_{2\ell+1}+o_{a.s.}(1)).
\end{align*}
This completes the proof.  $\blacksquare$

\bigskip
\bigskip

\noindent{\large\bf A.4\quad Proof of Theorem 2}

\medskip

\noindent Let $\bM(f_{t_k})=E[\bY_k\bY_k^T|f_{t_k}]$. Then the
matrix function $\bM(y)$ can be expanded around a fixed point $x$ as
$$\bM(y)=\bA_0+\bA_1(y-x)+\bA_2(y-x)^2+\bA_3(y-x)^3+\cdots,$$
where $\bA_0,\ \bA_1, \cdots$ are all matrices. To prove the
asymptotic property of the state-domain estimator, let us decompose
it as
\begin{align}\label{pf17}
\nonumber\hat{\bSig}_{S,t}(x)-\bM(x)& =\sum_{k=0}^{N-1} w_k(x)
\left(\bM(f_{t_k})-\bM(x)\right)
 + \sum_{k=0}^{N-1} w_k(x) \big(\bY_k\bY_k^T-\bM(f_{t_k})\big)\\
 &\equiv \bf{b}+\bf{t}.
\end{align}

First, we establish the asymptotic behavior of the bias term
$\bf{b}$. Applying Taylor's expansion and Proposition \ref{P2}
results in
\begin{align*}
{\bf b}&=\sum_{k=0}^{N-1} w_k(x)
\left(\bM(f_{t_k})-\bM(x)\right)\\
&=\sum_{k=0}^{N-1} w_k(x) \bA_1(f_{t_k}-x)+\sum_{k=0}^{N-1} w_k(x)
\bA_2(f_{t_k}-x)^2+o_{a.s.}(h^3)\\
&=h^2\mu_2\bA_2+o_{a.s.}(h^2).
\end{align*}
Since we have the following decomposition
$$\hat{\bSig}_{S,t}(x)-\bSig(x)=\big(\hat{\bSig}_{S,t}(x)-\bM(x)\big)+\big(\bM(x)-\bSig(x)\big)=
[\bf{b}+ \big(\bM(x)-\bSig(x)\big)]+\bf{t}, $$ and
$\bM(x)-\bSig(x)=o_p(\Delta)$, the asymptotic bias of the
state-domain estimator is
\begin{align}\label{pf18}
{\bf b}+
\big(\bM(x)-\bSig(x)\big)=\frac{1}{2}h^2\mu_2\bSig''(x)+o_{a.s.}(h^2)+o_p(\Delta).
\end{align}

Then, let us consider the variance term $\bf t$. Since $\bf t$ is a
matrix, we first vectorize it and then consider the asymptotic
normality of its linear combination, i.e. we look at the statistic
\begin{align*}
\tilde{\bf t}= \ba^T \text{vech}\Big(\sum_{k=0}^{N-1}
w_k(x)\big(\bY_k\bY_k^T-\bM(f_{t_k})\big)\Big),
\end{align*}
where $\ba$ is a constant vector. By Proposition \ref{P2},
\begin{align}\label{pf19}
\tilde{\bf t}&= \frac{1}{p(x)N}\sum_{k=0}^{N-1} K_h(f_{t_k}-x) a^T
\text{vech}\big(\bY_k\bY_k^T-\bM(f_{t_k})\big)\{1+o_{a.s.}(1)\}\\
\nonumber &\equiv\mathcal{A}_N\{1+o_{a.s.}(1)\}.
\end{align}
Therefore, we only need to show the asymptotic normality of
$\mathcal{A}_N$.  To this end, first let
$\vartheta_{N,k}=K_h(f_{t_k}-x) \ba^T
\text{vech}\big(\bY_k\bY_k^T-\bM(f_{t_k})\big)$. Then
$\mathcal{A}_N=\frac{1}{p(x)N}\sum_{k=0}^{N-1}\vartheta_{N,k}$.
Straightforward calculations give
\begin{align}\label{009}
\var(\vartheta_{N,k})&=E\left(K_h(f_{t_k}-x) \ba^T
\text{vech}\left(\bY_k\bY_k^T-\bM(f_{t_k})\right)\right)^2\\
\nonumber&=E\big\{K_h^2(f_{t_k}-x)E\Big[\big( \ba^T
\text{vech}(\bY_k\bY_k^T-\bM(f_{t_k}))\big)^2|f_{t_k}\Big]\big\}\\
\nonumber&=2E\left\{K_h^2(f_{t_k}-x)\left( \ba^T
P_D\bSig(f_{t_k})\otimes\bSig(f_{t_k})P_D^T\ba \right)\right\}\\
\nonumber&=2h^{-1}\nu_0p(x)\ba^T
P_D\bSig(x)\otimes\bSig(x)P_D^T\ba(1+o(1)),
\end{align}
where the last step follows from Taylor's expansion.

Note that 
$\bY_{t_\ell}$
only depends on the sample path of $f_t$ over time interval
$[t_{\ell}, \ t_{\ell+1}]$. Thus by conditioning on
$\mathcal{F}_{t_\ell}$, we obtain
\begin{align}\label{007}
\cov(\vartheta_{N,1},\vartheta_{N,\ell+1})=E\big[\vartheta_{N,1}K_h(f_{t_\ell}-x)
E\big(\ba^T
\text{vech}\big(\bY_{\ell}\bY_{\ell}^T-\bM(f_{t_\ell})\big)\big|\mathcal{F}_{t_\ell}\big)\big]=0,\
\ell \geq 1.
\end{align}
Combining (\ref{009}) and (\ref{007}) entails
\begin{align*}
\var(\mathcal{A}_N)=\frac{2\nu_0}{Nhp(x)}\ba^T
P_D\bSig(x)\otimes\bSig(x)P_D^T\ba(1+o(1)).
\end{align*}

Since a stationary Markov process satisfying the $G_2$ condition of
Rosenblatt (1970) is $\rho$-mixing, we can use ``big-block and
small-block" arguments similar to those used by Fan and Yao (2003,
Theorem 2.22, p.77) to prove the asymptotic normality of
$\mathcal{A}_N$. The lengthy details are omitted here. Thus,
\begin{align*}
\sqrt{Nh}\mathcal{A}_N\toD \mathcal{N}(0,2\nu_0p(x)^{-1}\ba^T
P_D\bSig(x)\otimes\bSig(x)P_D^T\ba).
\end{align*}
This together with (\ref{pf18}) and (\ref{pf19}) implies the
asymptotic normality of the state-domain estimator, i.e.
\begin{align*}
\sqrt{Nh}\ba^T\text{vech}\Big(\hat{\bSig}_{S,t}(x)-\bSig(x)-\frac{1}{2}h^2\mu_2\bSig''(x)\Big)\toD
\mathcal{N}(0,2\nu_0p(x)^{-1}\ba^T \Lambda(x)\ba),
\end{align*}
where $\ba$ is an arbitrary constant vector. This completes the
proof.  $\blacksquare$

\bigskip
\bigskip

\noindent{\large\bf A.5\quad Proof of Theorem 3}

\medskip

\noindent We only need to show the asymptotic normality of the
linear combination
$$\sqrt{Nh}\ \ba^T\vech\left(\widehat{\bSig}_{S,t}-\bSig(x)-\frac{1}{2}h^2\mu_2\bSig''(x)\right)
+ \sqrt{n}\
\bc^T\vech\left(\widehat{\bSig}^2_{T,t}-\bSig(x)\right),$$ where
$\ba^T$ and $\bc^T$ are two constant vectors. 
This is equivalent to showing the joint asymptotic normality of
$\sqrt{Nh}\ba^T\vech\big(\widehat{\bSig}_{S,t}-\bSig(x)-\frac{1}{2}h^2\mu_2\bSig''(x)\big)$
and $\sqrt{n}\bc^T\vech\big(\widehat{\bSig}^2_{T,t}\big)$.
From the proof of Theorem \ref{T2}, we have
\[
\ba^T\vech\big(\widehat{\bSig}_{S,t}-\bSig(x)-\frac{1}{2}h^2\mu_2\bSig''(x)\big)=
\ba^T \bt+
o_p(1)=\tilde{\bt}+o_p(1)=\mathcal{A}_N\{1+o_{a.s.}(1)\}+o_p(1),
\]
where $\bt$, $\tilde{\bt}$ and $\mathcal{A}_N$ are all defined in
the proof of Theorem 2. Therefore, we need only to consider about
the asymptotic normality of $\sqrt{Nh}\mathcal{A}_N$ and
$\sqrt{n}\bc^T\vech\big(\widehat{\bSig}^2_{T,t}\big)$.

We truncate $\mathcal{A}_N$ by defining
\[ \mathcal{A}_N^t=\frac{1}{p(x)N}\sum_{k=0}^{N-a_N}\vartheta_{N,k},
\]
where $a_N$ is an integer depending only on $N$ and satisfying
$a_N/N\rightarrow0$ and $a_N\Delta\rightarrow\infty$. We are going
to show that:
\begin{itemize}
\item[(i)] $\mathcal{A}_N^t$ and $\sqrt{n}\bc^T\vech\big(\widehat{\bSig}^2_{T,t}\big)$  are asymptotically
independent;
\item[(ii)] $\mathcal{A}_N-\mathcal{A}_N^t$ is asymptotically negligible.
\end{itemize}
We first prove (i). Since a stationary Markov process satisfying the
$G_2$ condition of Rosenblatt (1970) is $\rho$-mixing with
exponentially decaying $\rho$-mixing coefficient $\rho_t(\cdot)$,
and the strong-mixing coefficient $\alpha(\ell)\leq \rho(\ell)$ for
any integer $\ell$, it follows that
\[
\big|E\exp\{i\xi(\mathcal{A}_N^t+\bc^T\vech\big(\widehat{\bSig}^2_{T,t}\big))\}
-E\exp\{i\xi(\mathcal{A}_N^t)\}E\exp\{\i\xi\bc^T\vech\big(\widehat{\bSig}^2_{T,t}\big)\}\big|\leq
32\alpha(a_N-n)\rightarrow0,
\]
for any $\xi\in \mathbb{R}$. This proves (i).

Now, we prove (ii). From the proof of Theorem \ref{T2} we know that
\[ \var(\vartheta_{N,k})=2h^{-1}\nu_0p(x)\ba^T
P_D\bSig(x)\otimes\bSig(x)P_D^T\ba(1+o(1)), \] and
$\cov(\vartheta_{N,1},\vartheta_{N,\ell+1})=0$, $\forall \ell\geq1$.
Therefore,
\begin{align*}
\var(\sqrt{Nh}[\mathcal{A}_N-\mathcal{A}_N^t])=\frac{2a_N}{p(x)N}\nu_0\ba^T
P_D\bSig(x)\otimes\bSig(x)P_D^T\ba(1+o(1))\rightarrow0.
\end{align*}
This along with $E[\vartheta_{N,k}]=0$ gives
\[
\sqrt{Nh}[\mathcal{A}_N-\mathcal{A}_N^t]\toP0,
\]
which completes the proof of (ii). Combining (i) and (ii) entails
that $\sqrt{Nh}\mathcal{A}_N$ and
$\sqrt{n}\bc^T\vech\big(\widehat{\bSig}^2_{T,t}\big)$ are
asymptotically independent. This together with Theorem \ref{T1} and
the asymptotical normality of $\sqrt{Nh}\mathcal{A}_N$ shown in the
proof of Theorem \ref{T2} completes the proof of Theorem \ref{T3}.
$\blacksquare$

\bigskip

\newpage

\noindent{\Large\bf FIGURE LEGENDS}

\bigskip

\noindent{\bf\large Figure 1}. Illustration of time- and
state-domain estimation. (a)  The yields of 1-year, 5-year, and
10-year treasury bills from 1962 to 2005. The vertical bar indicates
localization in time, and the horizontal bar represents localization
in state of the 5-year treasury bill process. (b) Illustration of
time-domain smoothing: 1-year yield differences are plotted against
10-year yield differences with the regression line superimposed. (c)
Illustration of the state-domain smoothing: 1-year yield differences
are plotted against 10-year yield differences for those periods with
the corresponding 5-year yields restricted to the interval $6.37\%
\pm .2\%$, indicated by the horizontal bar in (a).

\bigskip

\noindent{\bf\large Figure 2}. Functions $A(\tau)$ (solid curve) and
$B(\tau)$ (dashed curve) for the parameters given in the simulation.

\bigskip

\noindent{\bf\large Figure 3}. (a) The averages of the entropy
losses over 500 simulations for the time-domain estimation (dotted
curve), state-domain estimation (dashed curve), and aggregated
method (solid curve). (b) The standard deviations of the entropy
losses over 500 simulations for time-domain estimation (dotted
curve), state-domain estimation (dashed curve), and the aggregated
method (solid curve). (c) and (d): The same as in (a) and (b) except
using the quadratic loss.

\bigskip

\noindent{\bf\large Figure 4}. (a) Box plots of the entropy losses
over 500 simulations for the time-domain estimator (left), the
aggregated method (middle), and the state-domain estimator (right).
(b) and (c): The same as in (a) except that the quadratic loss and
PE are used, respectively. (d) The ratios of the averages of the
quadratic losses over 150 out-sample forecastings using the
time-domain and state-domain estimators to those based on the
aggregated estimator ($x$-axis) are plotted against the ratios of
the PEs based on the time-domain and state-domain estimators to
those based on the aggregated estimator ($y$-axis).

\bigskip

\noindent{\bf\large Figure 5}. Correlation of the time-domain
estimator and state-domain estimator for the volatility of an
equally weighted portfolio. The dashed curves are for the 95\%
confidence intervals. The straight lines are acceptance regions for
testing the null hypothesis that the correlation is zero at
significance level 5\%.


\newpage
\begin{center}
{\normalsize REFERENCES}
\end{center}
\begin{singlespace}
\small
\begin{itemize}
\item[] A\"it-Sahalia, Y. (1996). ``Nonparametric Pricing
of Interest Rate Derivative Securities." {\em Econometrica} 64,
527--560.

\item[] A\"{\i}t-sahalia, Y., and P. Mykland. (2003). ``The Effects
of Random and Discrete Sampling When Estimating Continuous-Time
Diffusions." {\em Econometrica} 71, 483--549.

\item[]------ (2004). ``Estimating Diffusions with Discretely and Possibly Randomly
Spaced Data: A General Theory." {\em Annals of Statistics} 32,
2186--2222.

\item[] Andersen, T. G., T. Bollerslev, and F. X.
Diebold. (2002). ``Parametric and Nonparametric Volatility
Measurement," in {\em Handbook of Financial Econometrics} (Y.
A\"it-Sahalia and L. P. Hansen, eds.).

\item[] Arapis, M., and J. Gao. (2004). ``Nonparametric Kernel
Estimation and Testing in Continuous-Time Financial Econometrics."
{\em Manuscript}.

\item[] Arfi, M. (1998). ``Non-Parametric Variance Estimation from Ergodic
Samples." {\em Scandinavian Journal of Statistics} 25, 225--234.

\item[] Bandi, F. M., and G. Moloche. (2004). ``On the Functional Estimation
of Multivariate Diffusion Processes." {\em Manuscript}.

\item[] Bandi, F. M., and T. Nguyen. (1999). ``Fully Nonparametric
Estimators for Diffusions: A Small Sample Analysis." Working Paper,
University of Chicago.

\item[] Bandi, F. M., and P. C. B. Phillips. (2002). ``Nonstationary
Continuous-Time Processes," in {\em Handbook of Financial
Econometrics} (Y. A\"it-Sahalia and L. P. Hansen, eds.).

\item[] ------ (2003). ``Fully Nonparametric Estimation of Scalar Diffusion Models." {\em
Econometrica} 71, 241--283.

\item[] Banon, G. (1978). ``Nonparametric Identification
for Diffusion Processes." {\em SIAM J. Control Optim.} 16, 380--395.

\item[] Bollerslev, T., R. F. Engle, and J. M.
Wolldridge. (1988). ``A Capital Asset Pricing Model with
Time-Varying Covariance." {\em Jour. of Political Economy} 96,
116--131.

\item[] Cai, Z., and Y. Hong. (2003). ``Nonparametric Methods
in Continuous-Time Finance: A Selective Review." In {\em Recent
Advances and Trends in Nonparametric Statistics} (M. G. Akritas and
D. M. Politis, eds.), 283--302.

\item[] Chen, S.X. and J. Gao. (2004). ``A Test for Model Specification
of Diffusion Processes." {\em Manuscript}.

\item[] Cheridito, P., D. Filipovi\'{c}, and R. L.
Kimmel. (2005). ``Market Price of Risk Specification for Affine
Models: Theory and Evidence." {\em Journal of Financial Economics},
Forthcoming.

\item[] Cox, J. C., J. E. Ingersoll, and S. A. Ross.
(1985). ``A Theory of the Term Structure of Interest Rates." {\em
Econometrica} 53, 385--467.


\item[] Dalalyan, A. S., and Y. A. Kutoyants. (2003). ``Asymptotically
Efficient Estimation of the Derivative of the Invariant Density."
{\em Statist. Inference Stochastic Process.} 6, 89--107.

\item[] Duffee, G. R. (2002). ``Term Premia and Interest
Rate Forecasts in Affine Models." \emph{Journal of Finance} 57,
405--443.

\item[] Duffie, D., and R. Kan. (1996). ``A Yield-Factor Model of Interest
Rates." {\em Math. Finance} 6, 379?-406.

\item[] Engle, R. F., V. K. Ng, and M. Rothschild.
(1990). ``Asset Pricing with a Factor ARCH Covariance Structure:
Empirical Estimates for Treasury Bills." {\em Journal of
Econometrics} 45, 213--237.

\item[] Fan, J. (2005). ``A Selective Overview of Nonparametric Methods in
Financial Econometrics (with discussion)." {\em Statistical
Science}, 20, 317--357.

\item[] Fan, J., Y. Fan, and J. Jiang. (2005). ``Dynamic
Integration of Time- and State-Domain Methods for Volatility
Estimation." {\em Manuscript}.


\item [] Fan, J., J. Jiang, C. Zhang, and Z. Zhou. (2003).
``Time-Dependent Diffusion Models for Term Structure Dynamics and
the Stock Price Volatility." {\em Statistica Sinica} 13, 965--992.


\item[] Fan, J. and Q. Yao (2003). {\em Nonlinear Time Series: Nonparametric
and Parametric Methods}. New York: Springer-Verlag.

\item[] Fan, J., and C. Zhang. (2003). ``A Re-examination
of Stanton's Diffusion Estimations with Applications to Financial
Model Validation." {\em J. Amer. Statist. Assoc.} 98, 118--134.

\item[] Foster, D. P., and D. B. Nelson. (1996). ``Continuous
Record Asymptotics for Rolling Sample Variance Estimators." {\em
Econometrica} 64, 139--174.


\item[] Gobet, E. (2002). ``LAN Property for Ergodic Diffusions with Discrete
Observations." {\em Ann. Inst. H. Poincar\'{e} Probab. Statist.} 38,
711--737.

\item[] Gobet, E., M. Hoffmann, and M. Reiss. (2004). ``Nonparametric
Estimation of Scalar Diffusions Based on Low Frequency Data Is
Ill-Posed." {\em Ann. Statist.} 32, 2223--2253.

\item[] Hall, P., and C. Heyde. (1980). {\em Martingale
Limit Theorem and Its Applications.} Academic Press.

\item[] Hansen, L. P., and Scheinkman, J. A. (1995). ``Back to the
Future: Generating Moment Implications for Continuous-Time Markov
processes." {\em Econometrica} 63, 767--804.

\item[] Hansen, L. P., J. A. Scheinkman, and N. Touzi. (1998). ``Spectral
Methods for Identifying Scalar Diffusions." {\em Journal of
Econometrics} 86, 1--32.

\item[] H\"ardle, W. , H. Herwartz, and V. Spokoiny.
(2002). ``Time Inhomogeneous Multiple Volatility Modelling." {\em
Jour. Fin. Econometrics} 1, 55-95.

\item[] Jacod, J. (1997). ``Nonparametric Kernel Estimation
of the Diffusion Coefficient of a Diffusion." Pr\'{e}publication N.
405 du Laboratoire de Probabilit\'{e}s de l'Universit\'{e} Paris VI.

\item[] Jiang, G. J., and J. Knight. (1997). ``A Nonparametric
Approach to the Estimation of Diffusion Processes, with an
Application to a Short-Term Interest Rate Model." {\em Econometric
Theory} 13, 615--645.

\item[] Kessler, M., and M. S{\o}rensen. (1999). ``Estimating
Equations Based on Eigenfunctions for a Discretely Observed
Diffusion Process." {\em Bernoulli} 5, 299--314.

\item[] Karatzas, I., and S. Shreve. (1991). {\em Brownian Motion
and Stochastic Calculus} (2nd ed.). New York: Springer-Verlag.

\item[] Ledoit, O. and M. Wolf. (2003). ``Improved Estimation of
the Covariance Matrix of Stock Returns with an Application to
Portfolio Selection." {\em Journal of Empirical Finance} 10,
603--621.


\item[]  Mercurio, D. and V. Spokoiny. (2004). ``Statistical Inference
for Time-Inhomogeneous Volatility Models." {\em The Annals of
Statistics} 32, 577--602.

\item[] Morgan, J. P. (1996). {\em RiskMetrics Technical Document}
(4th ed.). New York.

\item[] Revuz, D., and M. Yor. (1999). {\em Continuous
Martingales and Brownian Motion}. Springer-Verlag.

\item[] Robinson, P. M. (1997). ``Large-sample Inference for
Nonparametric Regression with Dependent Errors." {\em Ann. Statist.}
25, 2054--2083.

\item[] Rosenblatt, M. (1970). ``Density Estimates and
Markov Sequences," in {\em Nonparametric Techniques in Statistical
Inference }(M. L. Puri, ed.), 199--213. Cambridge University Press.


\item[] Stanton, R. (1997). ``A Nonparametric Models of Term
Structure Dynamics and the Market Price of Interest Rate Risk." {\em
Journal of Finance} 52, 1973--2002.

\item[] Vasicek, O. A. (1977). ``An Equilibrium Characterisation
of the Term Structure." \textit{Journal of Financial Economics} 5,
177--188.

\end{itemize}
\end{singlespace}

\newpage

\noindent{\Large\bf FOOTNOTE}

\bigskip

\noindent{\bf\large Footnote 1}. By ``stationarity" we do not mean
that the process is strongly stationary, but has some structural
invariability over time. For example, the conditional moment
functions do not vary over time.

\smallskip

\noindent{\bf\large Footnote 2}. Ledoit and  Wolf (2003) introduce a
shrinkage estimator by combining the sample covariance estimator
with that derived from the CAPM.  Their procedure intends to improve
estimated covariance matrix by pulling the sample covariance towards
the estimate based on the CAPM.  Their basic assumption is that the
return vectors are i.i.d. across time.  This usually holds
approximately when the data are localized in time.  In this sense,
their estimator can be regarded as a time-domain estimator.

\smallskip

\noindent{\bf\large Footnote 3}. We prove in Section 4 that
$\widehat{\bSig}_{S,t}$ and $\widehat{\bSig}_{T,t}$ are
asymptotically independent, and thus they are close to be
independent in finite sample. In the following, by ``nearly
independent'' and ``almost uncorrelated'', we mean the same.

\smallskip

\noindent{\bf\large Footnote 4}. In practice, one can take the
yields process with median term of maturity as the driving factor,
as this bond is highly correlated to both short-term and long-term
bonds.

\smallskip

\noindent{\bf\large Footnote 5}. The kernel function is a
probability density, and the bandwidth is its associated scale
parameter.  Both of them are used to localize the linear regression
around the given point $x_0$.  The commonly used kernel functions
are the Gaussian density and the Epanechnikov kernel $K(x) = 0.75
(1-x^2)_+$.

\smallskip

\noindent{\bf\large Footnote 6}. The stationarity condition of $f_t$
in Assumption 3 can be weakened to Harris recurrence. See Bandi and
Moloche (2004) for asymptotic normality of local constant estimator
under recurrence assumption.

\smallskip

\noindent{\bf\large Footnote 7}. The optimal choice of weight is
proportional to the effective number of data points used for the
state-domain and time-domain smoothing. It always outperforms the
choice with $\omega_t=1$ (state-domain estimator) or $\omega_t = 0$
(time-domain estimator).

\smallskip

\noindent{\bf\large Footnote 8}. The choice comes from the
recommendation of the RiskMetrics of J.P. Morgan.  The parameter
$\lambda$ can also be chosen automatically by data by using the
prediction error as in Fan, Jiang, Zhang and Zhou (2003).  Since we
compare the relative performance between the time-domain estimator
and the aggregated estimator, we opt for this simple choice.  The
results do not expect to change much when a data-driven technique is
used.

\smallskip

\noindent{\bf\large Footnote 9}. Here we add normal noise to make
the model more realistic. Our method performs even better without
noise. Since the noise vectors are i.i.d. across time and the
standard deviations are small, adding them to the original time
series does not change the whole structure. Hence, our theory can
carry through under contamination.

\smallskip

\noindent{\bf\large Footnote 10}. With $\lambda=0.94$, the last data
point used in the time domain has an extra weight $0.94^{104}\approx
0.0016$, which is very small. Hence, we essentially include all the
effective data points.

\smallskip

\noindent{\bf\large Footnote 11}. Europe used several common
currencies prior to the introduction of the Euro. The European
Currency Unit (ECU) was used from January 1, 1979 to January 1,
1999, when the Euro replaced the European Currency Unit at par.

\vspace{1 in}

\noindent{\Large\bf TABLE}

\bigskip

\begin{center}
\centering{\small {\bf\large Table 1}\\
\textsc{\small APEs of Bond Yields, Exchange Rates and Simulations}}\\
\vspace{0.1 in}
\begin{tabular}{c|ccc}
  \hline
 & \text{Time} & \text{State} & \text{Aggregated} \\
  \hline
\text{Bonds}\\
 $k=0$   &$3.837\times10^{-3}$  & $3.767\times10^{-3}$  & $3.756\times10^{-3}$\\
  $k=1$   &$1.643\times10^{-3}$  & $1.557\times10^{-3}$  & $1.555\times10^{-3}$\\
 $k=2$   &$1.013\times10^{-3}$  & $1.011\times10^{-3}$  & $9.933\times10^{-4}$\\
  \hline
  \text{Currencies}\\
  $k=0$   &$4.795\times10^{-3}$  & $4.913\times10^{-3}$  & $4.755\times10^{-3}$\\
 $k=1$   &$1.681\times10^{-3}$  & $1.855\times10^{-3}$  & $1.652\times10^{-3}$\\
 $k=2$   &$8.979\times10^{-4}$  & $1.184\times10^{-3}$  & $8.937\times10^{-4}$\\
  \hline
  \text{Simulations ($k=0$)} &$1.850\times10^{-2}$  & $1.846\times10^{-2}$&
$1.825\times10^{-2}$\\
\hline
   \end{tabular}
\end{center}

\newpage

\begin{singlespace}
\begin{figure} \centering
\begin{center}%
\begin{tabular}
[c]{c}%
{\includegraphics[ trim=0.000000in 0.000000in 0.000000in
-0.189256in, height=3.5in, width=4.5in
]%
{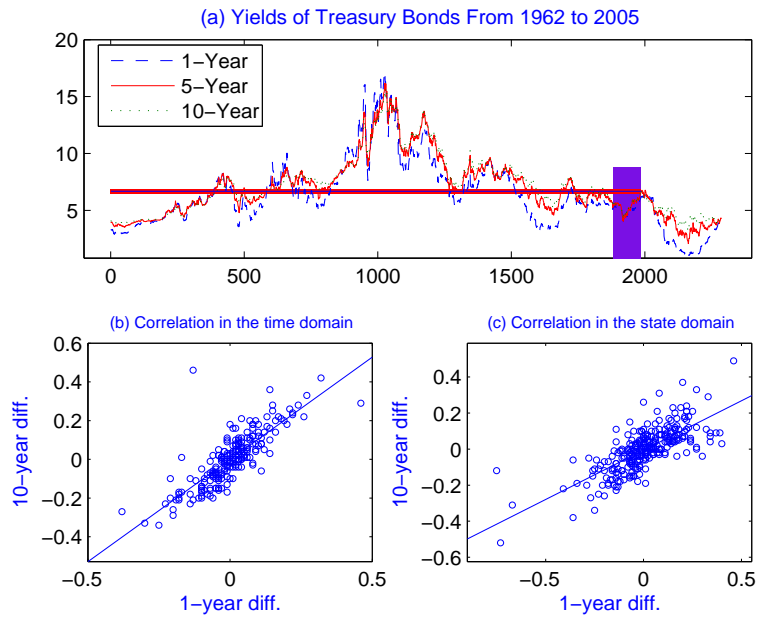}%
}%
\end{tabular}%
\caption{\scriptsize  Illustration of time- and state-domain
estimation. (a)  The yields of 1-year, 5-year and 10-year treasury
bills from 1962 to 2005. The vertical bar indicates localization in
time, and the horizontal bar represents localization in the state of
the 5-year treasury bill process. (b) Illustration of time-domain
smoothing: 1-year yield differences are plotted against 10-year
yield differences with the regression line superimposed. (c)
Illustration of the state-domain smoothing: 1-year yield differences
are plotted against 10-year yield differences for those periods with
the corresponding 5-year yields restricted to the interval $6.37\%
\pm
.2\%$, indicated by the horizontal bar in (a).}%
\end{center}%
\end{figure}%
\end{singlespace}

\begin{singlespace}
\begin{figure} \centering
\begin{center}%
\begin{tabular}
[c]{c}%
{\includegraphics[ trim=0.000000in 0.000000in 0.000000in
-0.189256in, height=3in, width=4in
]%
{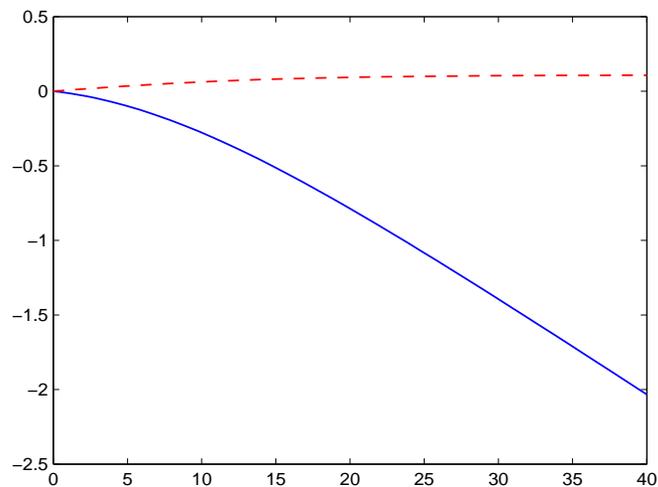}%
}%
\end{tabular}%
\caption{\scriptsize Functions $A(\tau)$ (solid curve) and
$B(\tau)$ (dashed curve) for the parameters given in the simulation.}%
\end{center}%
\end{figure}%
\end{singlespace}

\newpage

\begin{singlespace}
\begin{figure} \centering
\begin{center}%
\begin{tabular}
[c]{cc}%
{\includegraphics[ trim=-0.000000in 0.000000in 0.000000in
-0.189256in, height=2.2in, width=2.6in
]%
{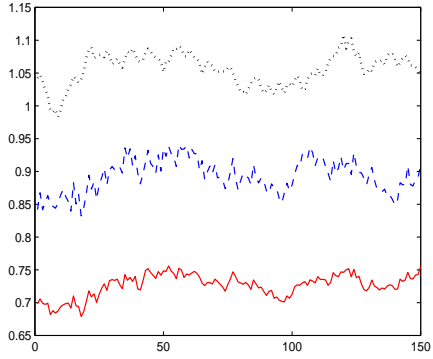}%
}%
& {\includegraphics[ trim=0.000000in 0.000000in 0.000000in
-0.189256in, height=2.2in, width=2.6in
]%
{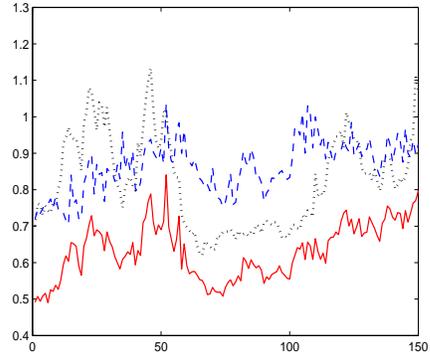}%
}%
\\
\scriptsize (a) & \scriptsize (b) \\
 {\includegraphics[ trim=0.000000in 0.000000in
0.000000in -0.189256in, height=2.2in, width=2.6in
]%
{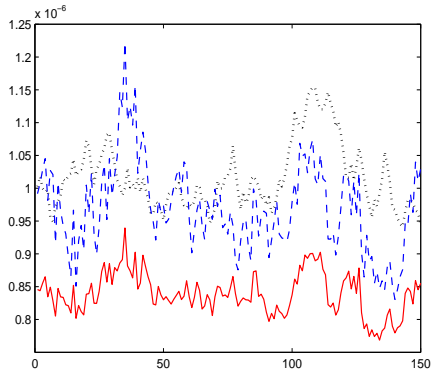}%
}%
& {\includegraphics[ trim=0.000000in 0.000000in 0.000000in
-0.189256in, height=2.2in, width=2.6in
]%
{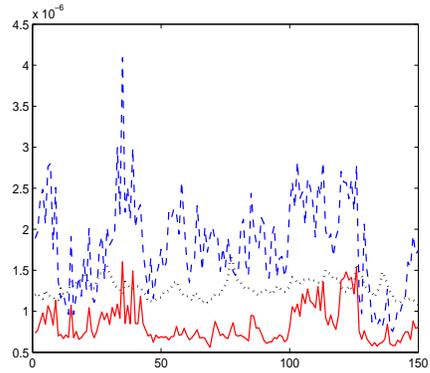}%
}%
\\ \scriptsize (c) & \scriptsize (d)
\end{tabular}%
\caption{\scriptsize (a) The averages of the entropy losses over 500
simulations for the time-domain estimation (dotted curve),
state-domain estimation (dashed curve) and aggregated method (solid
curve). (b) The standard deviations of the entropy losses over 500
simulations for the time-domain estimation (dotted curve),
state-domain estimation (dashed curve) and aggregated method (solid
curve). (c) and (d): The same as in (a) and (b) except using the
quadratic loss.}%
\end{center}%
\end{figure}%
\end{singlespace}

\newpage

\begin{singlespace}
\begin{figure} \centering
\begin{center}%
\begin{tabular}
[c]{cc}%
{\includegraphics[ trim=0.000000in 0.000000in 0.000000in
-0.189256in, height=2.2in, width=2.6in
]%
{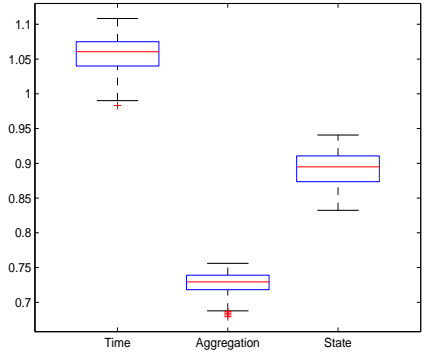}%
}%
& {\includegraphics[ trim=0.000000in 0.000000in 0.000000in
-0.189256in, height=2.2in, width=2.6in
]%
{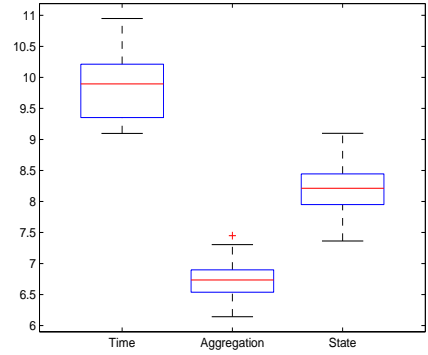}%
}%
\\
\scriptsize (a) & \scriptsize (b) \\
 {\includegraphics[ trim=0.000000in 0.000000in
0.000000in -0.189256in, height=2.2in, width=2.6in
]%
{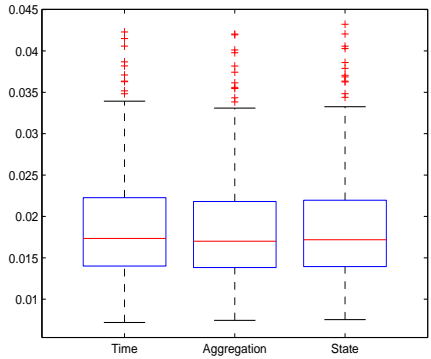}%
}%
& {\includegraphics[ trim=0.000000in 0.000000in 0.000000in
-0.189256in, height=2.2in, width=2.6in
]%
{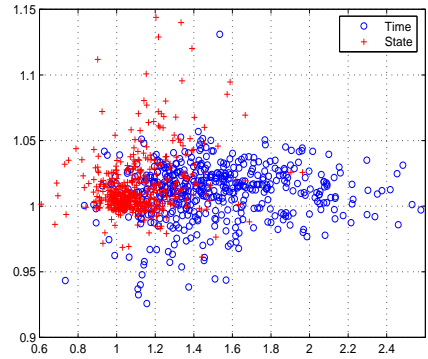}%
}%
\\ \scriptsize (c) & \scriptsize (d)
\end{tabular}%
\caption{\scriptsize (a) Box plots of the entropy losses over 500
simulations for the time-domain estimator (left), the aggregated
method (middle), and the state-domain estimator (right). (b) and
(c): The same as in (a) except that the quadratic loss and PE are
used, respectively. (d) The ratios of the averages of the quadratic
losses over 150 out-sample forecastings using the time-domain and
state-domain estimators to those based on the aggregated estimator
($x$-axis) are plotted against the ratios of the PEs based on the
time-domain and state-domain estimators to those based on the
aggregated estimator ($y$-axis).}
\end{center}%
\end{figure}%
\end{singlespace}

\newpage

\begin{singlespace}
\begin{figure} \centering
\begin{center}%
\includegraphics[ trim=0.000000in 0.000000in 0.000000in
-0.189256in, height=2.2in, width=2.6in
]%
{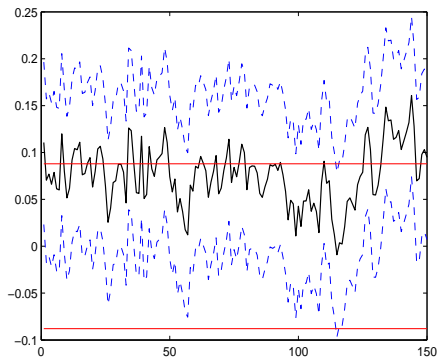}%
\caption{\scriptsize Correlation of the time-domain estimator and
state-domain estimator for the volatility of an equally weighted
portfolio. The dashed curves are for the 95\% confidence intervals.
The straight lines are acceptance regions for testing the null
hypothesis
that the correlation is zero at significance level 5\%.}%
\end{center}
\end{figure}
\end{singlespace}

\end{document}